\documentclass{amsart}
\usepackage{amsmath,amsthm,amssymb,amsfonts,amscd}
\usepackage{mathrsfs}
\usepackage{bbding}
\usepackage{graphicx,latexsym}
\usepackage{hyperref}
\usepackage{geometry}
\usepackage{color}
\usepackage{xcolor}
 \usepackage{floatrow}

\hypersetup{
    colorlinks,
    linkcolor={red!50!black},
    citecolor={blue!50!black},
    urlcolor={blue!80!black}
}

\numberwithin{equation}{section}

\theoremstyle{plain}
\newtheorem{theorem}{Theorem}
\newtheorem*{theorem*}{Theorem}
\newtheorem*{conjecture*}{Conjecture}
\newtheorem{lemma}{Lemma}
\newtheorem{corollary}{Corollary}
\newtheorem{proposition}{Proposition}

\theoremstyle{definition}

\newtheorem{conjecture}[equation]{Conjecture}

\theoremstyle{remark}
\newtheorem{remark}{Remark}

\renewcommand{\mod}{\operatorname{mod}\,}

\newcommand{\fE}{\mathfrak{E}}

\newcommand{\mat}[1]{\left(\begin{matrix} #1 \end{matrix} \right)}  

\renewcommand{\d}{\delta}

\newcommand{\ve}{\varepsilon}

\begin{document}

\title[Thin Bianchi groups]
        {An asymptotic local-global theorem on heights of some Kleinian group orbits }
\author{Xuanxuan Xiao}
\address{Macao Centre for Mathematical Sciences \\ Macau University of Science and Technology \\ Macau}
\email{xxxiao@must.edu.mo}
\author{Xin Zhang}
\address{Department of Mathematics \\ The University of Hong Kong}
\email{xzhang@maths.hku.hk}

\date{\today}

\begin{abstract} We use circle method prove an asymptotic local-global theorem on the heights of point orbits of thin subgroups of Bianchi groups in $\mathbb H^3$. 
\end{abstract}

\keywords{thin group, Bianchi groups, circle method, local-global principle}
\subjclass[2020]{11N25, 11N45, 11P21, 11P55, 11S23}
\thanks{The first author is supported The Science and Technology Development Fund, Macau SAR (File no. 0066/2020/A2).  The second author is partially supported by GRF grant 17317222 and NSFC grant 12001457.}

\maketitle

\section{Introduction}

In the papers \cite{BGP1,BGP2}, Bourgain, Gamburd and Sarnak initiated the study of the following general question: Let $\Gamma < \text{GL}_n(\mathbb Z)$, $v\in \mathbb Z^n, f\in \mathbb Z[x_1,\cdots,x_n]$ an integral polynomial of $n$ variables. What integers can be written as $f(\gamma\cdot v)$? This program includes Apollonian circle packing curvature problem and Zaremba's conjecture as well as other interesting problems as special cases. Using sieve method, Salehi-Golsefidy and Sarnak obtained a general theorem which says that under some fairly weak assumption of $\Gamma$, the representation set $f(\Gamma\cdot v)$ contains infinitely many pseudo-primes. Obtaining finer information beyond what sieve gives is often of a technical and ad-hoc endeavor. Remarkably, Bourgain-Kontorovich used circle method to prove a density one theorem for both the curvature problem and Zaremba's conjecture \cite{BK1,BK2}. Some other relevant works are \cite{FSZ,KZ23, ZHA,ZHA1}.
\par

In this paper we investigate a new incidence of this program. Let $$\mathbb H^3:=\{x+y{\bf i}+z{\bf j}:x,y\in \mathbb R, z\in\mathbb R_+\}$$
be the standard upper half space model for the hyperbolic 3-space. $\text{SL}_2(\mathbb C)$ acts on $\mathbb H^3$ as orientation preserving isometries by linear fractional transformations: 
$$\mat{a&b\\c&d}h=(ah+b)(ch+d)^{-1}, h\in \mathbb H^3.$$

Let $\Gamma$ be a subgroup of the Bianchi group $\text{SL}_2(\mathbb Z[\sqrt{-\mathfrak D}])$, where $\mathfrak D$ is a squarefree number. We fix an element in $\mathbb H^3$, say $\bf j$, and we look at the orbit $\Gamma\cdot {\bf j}$. The height, or the Euclidean distance of a point to the boundary of plane $\mathbb C$, gives one natural way to measure the ``closeness" of this point to infinity. It turns out the inverse of heights of this point orbit are all integers.  This follows from the explicit formula that if we write $\gamma\in \text{SL}_2(\mathbb C)$, then $$\gamma{\bf j}=\frac{ac+b\bar{d}}{|c|^2+|d|^2}+\frac{1}{|c|^2+|d|^2}{\bf j}.$$ 
So the inverse height $H(\gamma):=|c|^2+|d|^2\in\mathbb Z$ if $\gamma=\mat{a&b\\c&d}\in \text{SL}_2(\mathbb Z[\sqrt{-\mathfrak D}])$. In this paper we study the integer set $H(\Gamma)=\{H(\gamma):\gamma\in\Gamma\}$.

If $\Gamma$ has critical exponent $\delta=\delta(\Gamma)>1$, then $H(\Gamma)\cap[-N, N]$ counting multiplicity has size $\asymp N^{\delta}$, so the growth rate of $H(\Gamma)$ exceeds that of integers. It then makes sense to conjecture a local-global principle: 

\begin{conjecture}[Local-Global Conjecture] \label{1205} Let $\Gamma$ be a geometrically finite subgroup of $\text{SL}_2(\mathbb Z[\sqrt{-\mathfrak D}])$ with critical exponent $\delta_\Gamma>1$, then $H(\Gamma)=\{H(\gamma):\gamma\in\Gamma\}$ and its admissible set $\mathcal{A}(\Gamma)$ differ by a finite set.  
\end{conjecture}
\begin{remark}
The strong approximation property of $\Gamma$ implies that there is a finite number $L$ that describes the non-Archimidean obstruction of $H(\Gamma)$, i.e., an integer $n\in H (\mod q)$ for an arbitrary modulus $q$ if and only if $n\in H(\mod L)$. We call an congruence classes mod $L$ \emph{admissible} if it contains at least one element in $H$. There is also an obvious obstruction at infinity: all numbers in $H(\Gamma)$ are positive. The set $\mathcal{A}(\Gamma)$ is then the positive numbers in all admissible congruence classes mod $L$. 
\end{remark}
\begin{remark}
Note that when $\mathfrak D=1$, Conjecture \ref{1205} is a refinement of Lagrange's four square theorem. 
\par
\end{remark}

The most general case of Conjecture \ref{1205} seems out of reach from the current techniques. To obtain a definite result, we make some further assumptions on $\Gamma$. The following is our main theorem: 

\begin{theorem}\label{mainthm} Let $\Gamma$ be a geometrically finite subgroup of $\text{SL}_2(\mathbb Z[{\sqrt{-\mathfrak D}}])$, containing a rank-2 hyperbolic subgroup. There exists absolute constant $1<\delta_0<2$ such that if $\Gamma$ has critical exponent $\delta(\Gamma)>\delta_0$, 
$$\frac{\#H(\Gamma)\cap [0, N]}{\# \mathcal A(\Gamma)\cap [0, N]} =1+O(N^{-\eta}),$$
where $\eta>0$ depends only on the group $\Gamma$. 
\end{theorem}
\begin{remark}
We assume $\Gamma$ is geometrically finite so that we can apply counting techniques developed in \cite{ BK1, LO12, VIN14}.
\end{remark}

\begin{remark}
We assume that $\Gamma$ contains a rank two parabolic subgroup $P$, say $$P=\left\{\left\{\mat{1&x+y{\sqrt{-\mathfrak D}}\\0&1}: x,y\in\mathbb Z\right\}\right\},$$ so that integers represented by the (shifted) binary quadratic form $Q_\gamma(x,y)=H\left(\gamma\cdot\mat{1&x+y{\sqrt{-\mathfrak D}}\\0&1}\right)$ belong to the representation set $H(\Gamma)$. Since a single binary quadratic form represents $\gg N/\sqrt{\log N}$ integers, we can hope that a combination of different $Q_\gamma$ will give a density-one subset. We achieve this using the orbital circle method and the general framework laid out in \cite{BK}.
\end{remark}

\begin{remark}
The constant $\delta_0$ is absolute and very close to $2$. It depends on the spectral gap which is absolute according to \cite{MA15}. To our surprise,  we have the extra dimension requirement as stated in Theorem \ref{mainthm}. By comparison, the method in \cite{BK1}, \cite{FSZ} and \cite{KZ23} only requires that $\delta>1$, which is automatically guaranteed by the existence of surface subgroups showing up in these works. One reason is all the binary forms arising in this paper have equal quadratic coefficients, so more restrictive than the quadratic forms produced by the surface groups in the aforementioned preceding works. 
\end{remark}

\begin{remark}
The groups satisfying the assumption for $\Gamma$ in Theorem \ref{mainthm} are abundant.  Starting with a lattice $\Gamma_0< \text{SL}_2(\sqrt{-\mathfrak D})$ with at least two inequivalent rank-2 cusps, one can follow the construction in the proof of Proposition 60 in \cite{KK23} to construct a family of groups $\Gamma_j, j\in \mathbb N$, each of which has a rank-2 parabolic subgroups.  Theorem 13 in \cite{KK23} then guarantees that $\delta(\Gamma_j)\rightarrow 2$ as $j\rightarrow \infty$, so $\Gamma_j$ satisfies the assumption for $\Gamma$ in Theorem \ref{mainthm} for sufficiently large $j$. \par
\end{remark}

\begin{remark} Remarkable counterexamples to the full local-global conjecture for the Apollonian packing curvature problem and generalized Zaremba's conjecture are given in \cite{HKRS,RS} using quadratic reciprocity. We have not detected a quadratic obstruction in our setting, so we state Conjecture \ref{1205} as it is. The works \cite{HKRS,RS} expose an intrinsic limitation of the circle method, as the conventional circle method does not detect this quadratic obstruction, which already gives 1/2 power loss for the admissible set. It seems that the power-saving results of the works \cite{BK1,BK2,FSZ,KZ23} and our paper are close to the limit that circle method can achieve. 

\end{remark}

\subsection*{Acknowledegements}

We thank Alex Kontorovich for discussions about this work.

\newpage
\section{Setup of the Problem and The Circle Method}

\subsection{Setup of the problem}

Let 
\begin{equation*}       
\gamma=\left(                 
\begin{array}{ccc}   
a_1+a_2\sqrt{-\mathfrak D} & b_1+b_2\sqrt{-\mathfrak D} \\  
c_1+c_2\sqrt{-\mathfrak D} &  d_1+d_2\sqrt{-\mathfrak D}\\  
\end{array}
\right) 
\in \Gamma.                
\end{equation*}
One can compute that for $x,y\in \mathbb Z$
$$H\left(\gamma\cdot\mat{1&x+y\sqrt{-\mathfrak D}\\0&1}\right)=A_{\gamma}x^2+B_{\gamma}y^2+C_{\gamma}x+D_{\gamma}y+E_{\gamma}:=Q_{\gamma}(x,y),$$
with
\begin{equation}\label{defQ}
\begin{aligned}
A_{\gamma}&=c_1^2+c_2^2{\mathfrak D}, \quad B_{\gamma}=(c_1^2+c_2^2{\mathfrak D}){\mathfrak D},\\
C_{\gamma}&=2(c_1d_1+c_2d_2{\mathfrak D}), \quad D_{\gamma}=2(c_1d_2-c_2d_1){\mathfrak D},\\
E_{\gamma}&=c_1^2+c_2^2{\mathfrak D}+d_1^2+d_2^2{\mathfrak D}.
\end{aligned}
\end{equation}

\medskip

Let $\parallel \cdot \parallel$ be the standard norm on $\text{SL}_2(\mathbb{Z}[\sqrt{-\mathfrak D}])$ given by
$$\parallel \gamma \parallel=\sqrt{a_1^2+a_2^2{\mathfrak D}+b_1^2+b_2^2{\mathfrak D}+c_1^2+c_2^2{\mathfrak D}+d_1^2+d_2^2{\mathfrak D}}.$$

Let $N>0$ be the main growing parameter.  We choose two parameters $T$ and $X$, such that $T=N^{\sigma}$ with $\sigma$ small ($1/2024$ for example), and $X=N^{\frac{1}{2}-\sigma}$, so that $N=T^2X^2$. 
Denote by
$$B_T=\left\{\gamma \in \Gamma:  \parallel \gamma \parallel <T,A_\gamma \geq \frac{T^2}{100}\right\}.$$
Theorem 1.5 of \cite{BKS} implies that $\#B_T\asymp T^{2\delta}$.

For $N< n \leq 2N$, define the counting function
\begin{equation*}
R(n)=\sum_{\gamma \in B_T}\sum_{x,y \in \mathbb{Z}}\Psi(\frac{x}{X})\Psi(\frac{y}{X})1\{Q_{\gamma}(x,y)=n\},
\end{equation*}
where $\Psi\in C(\mathbb R)$ is a smooth function satisfying
\begin{itemize}
\item $\Psi\geq 0$,
\item Supp$(\Psi)\subset [0.5,2.5],$ 
\item $\Psi=1$ on $[1,2]$.
\end{itemize}
Thus, $R(n)>0$ implies that $n$ is represented. Note that $|A_{\gamma}|,|B_{\gamma}|\asymp T^2$, and 
\begin{align}\label{TXN}
|Q_{\gamma}(x,y)|\asymp T^2X^2=N.
\end{align}
Our total input is 
\begin{equation*}
\sum_{n\in \mathbb{Z}} R(n)\sim T^{2\delta}X^2\sim T^{2\delta-2}N.
\end{equation*}
We expect the total input is roughly equidistributed among admissible integers, i.e., for $n\sim N$ admissible, we expect 
$R(n)\sim T^{2\delta-2}$.

\subsection{The circle method}

We study the representation function $R(n)$ via its Fourier transformation
\begin{equation}\label{defRhat}
\widehat{R}(\theta)=\sum_{\gamma\in B_T}\sum_{x,y \in \mathbb{Z}}\Psi(\frac{x}{X})\Psi(\frac{y}{X})e(Q_{\gamma}(x,y)\theta).
\end{equation}
Then
\begin{equation*}
R(n)=\int_0^1\widehat{R}(\theta)e(-n\theta)d\theta.
\end{equation*}

We split the above integral into a sum of integrals over ``major'' and ``minor'' arcs
\begin{equation*}
R(n)=\int_0^1\widehat{R}(\theta)e(-n\theta)d\theta
=\int_{\mathcal{M}}+\int_{\mathfrak{m}}=: \mathcal{M}_N(n)+\mathcal{E}_N(n).
\end{equation*}

To define the major arc, we need the following wedge function 
$$\mathfrak{t}(x):=\min (1+x,1-x)^+,$$
 whose Fourier transform is 
$$\widehat{\mathfrak{t}}(y)=\left(\frac{\sin(\pi y)}{\pi y}\right)^2.$$
Let $Q_0$ and $K_0$ be a small power of $T$ chosen later (see \eqref{KQ}). Denote  
\begin{equation}\label{T}
\mathfrak{T}(\theta):=\sum_{q\leq Q_0}{\sum_{r (\mod q)}}'\sum_{m\in \mathbb{Z}}\mathfrak{t}\left(\frac{N}{K_0}(\theta+m-\frac{r}{q})\right).
\end{equation}
Then the function $\mathfrak{T}$ captures small neighborhoods of rationals with denominators $\leq Q_0$.

The main term is then defined as 
\begin{equation}\label{Major}
\mathcal{M}_N(n):=\int_0^1\mathfrak{T}(\theta)\widehat{R}(\theta)e(-n\theta)d\theta,
\end{equation}
and the minor term is defined as

\begin{equation*}
\mathcal{E}_N(n):=\int_0^1(1-\mathfrak{T}(\theta))\widehat{R}(\theta)e(-n\theta)d\theta.
\end{equation*}

In Section \ref{major} we show that for an admissible $n\asymp N$, we have $\mathcal{M}_N(n)\gg T^{2\delta-2}/\log N$ (Proposition \ref{main arcs}). In Section \ref{minor}, an ${\mathcal L}^2$ estimate for $\mathcal E_N$ (Propositition \ref{minor arcs}) is obtained, which implies most $\mathcal E_N(n)$ are small compared to $\mathcal M_N(n)$, leading to Theorem \ref{mainthm}.

\section{Some Lemmas}
In this section, we are going to list some lemmas which will be useful for the future use.

Let 
\begin{align*}
\Gamma(q)
=\left\{\left(                 
\begin{array}{ccc}   
a & b \\  
c &  d\\  
\end{array}
\right) 
\in \Gamma\, \Big| 
\left(                 
\begin{array}{ccc}   
a & b \\  
c &  d\\  
\end{array}
\right) 
\equiv I \  (\mod q)
\right\}
\end{align*}
be the principal congruence subgroup of $\Gamma$ of level $q$.\par
The first lemma is due to Bourgain-Kontorovich \cite{BK1} which gives a quantitative equidistribution of $\Gamma$ in each ``admissible cosets'' of $\Gamma(q)$: 
\begin{lemma}\label{lemmaBKS}
	Let $1<K<T^{1/10}$, $|\beta|<K/N$ and $x,y\asymp X$. Then for any $\gamma_0 \in \Gamma$ and $q\geq 1$, we have
	\begin{equation*}
	\sum_{\substack{\gamma \in B_T\\ \gamma\equiv \gamma_0 (\Gamma(q))}}e(\beta Q_{\gamma}(x,y))=\frac{1}{[\Gamma:\Gamma(q)]}\sum_{\gamma \in B_T}e\left(\beta Q_\gamma(x,y)\right)+O(T^{2\Theta_1}K),
	\end{equation*}
	where $\Theta_1<\delta$ depends only on the spectral gap for $\Gamma$ and the implied constant does not depend on $q,\gamma_0$, $x$ or $y$.
\end{lemma}

For $z\in \mathbb{Q}_p$, we define
$$\deg_{p^\nu}(z)=\max_{-\infty<k\leq \nu}\{k:p^{-k}z\in \mathbb{Z}_p\}.$$

The following two lemmas are direct corollaries of Gauss sums.
\begin{lemma}\label{gauss sum}
	Let $p$ be an odd prime and $\nu>0$ be an integer. For $a,b\in \mathbb{Z}$, we have
	\begin{align*}
	&\sum_{x\  (\mod p^{\nu})}e_{p^\nu}(ax^2+bx)\\
   =& \begin{cases}
   p^\nu \cdot {\bf 1}\{p^\nu \mid b\},& \text{if}\  p^\nu \mid a,\\
   p^{\nu/2}(p^\nu,a)^{1/2}i^{\epsilon(p^\nu/(p^\nu,a))}\left(\frac{\frac{a}{(p^{\nu-1},a)}}{p}\right)e_{p^\nu}\left(-\frac{b^2}{4a}\right)\cdot {\bf 1}\{\deg_{p^\nu}(b)\geq \deg_{p^\nu}(a)\},& \text{if} \  p^\nu \nmid a.
   \end{cases}
	\end{align*}
Here $\epsilon(n)=0$ if $n\equiv 1 \ (\mod 4)$ and $\epsilon(n)=1$ if $n\equiv 3 \ (\mod 4)$. 
\end{lemma}
In the above lemma, we use the notation that
$$e_{p^\nu}\left(-\frac{b^2}{4a}\right)=e\left(\frac{b'^2\overline{4a'}}{p^{\nu'}}\right)$$
with 
$a'=\frac{a}{\gcd(a,p^\nu)}$, $b'=\frac{b}{\gcd(a,p^\nu)}$ and $\overline{4a'} \cdot 4a' \equiv 1 (\mod p^{\nu'})$. When $p^\nu \mid a$, $e_{p^\nu}\left(-\frac{b^2}{4a}\right)=1$ since $p^{\nu'}=1$. Thus, it is easy to see that the case $p^\nu \mid a$ can be covered by the case $p^\nu \nmid a$.

\begin{lemma}\label{gauss sumeven}
	Let $q=2^\nu$. For $a,b\in \mathbb{Z}$, we have
	\begin{align*}
	&\sum_{x\  (\mod q)}e_{q}(ax^2+bx)\ll q\cdot {\bf 1}\{2^{\nu-1}\mid a\}.
	\end{align*}
\end{lemma}

We will deal with exponential sums of the form
\begin{align*}
S(q,A,B,C,D)=\sum_{x,y(\mod q)}e_q(Ax^2+By^2+Cx+Dy),
\end{align*} 
where $q,A,B,C,D$ are positive integers.
\begin{lemma}\label{ESQF}
	Let $p$ be a prime and $k_g=\deg_{p^\nu}(\gcd(A,B))$. 
	 If $p$ is odd and $k_g=\nu$, then
	$$S(p^\nu,A,B,C,D)=p^{2\nu}{\bf 1}\{p^\nu\mid \{C,D\}\}.$$
	If $k_g<\nu$, then
	\begin{equation*}
	\begin{aligned}
	S(p^\nu,A,B,C,D)=&p^\nu(p^\nu,A)^{1/2}(p^\nu,B)^{1/2}i^{\epsilon(p^\nu/(p^\nu,A))}i^{\epsilon(p^\nu/(p^\nu,B))}\\
	&\times \left(\frac{\frac{A}{(p^{\nu-1},A)}}{p}\right)\left(\frac{\frac{B}{(p^{\nu-1},B)}}{p}\right)\chi(p^\nu,A,B,C,D)e_{p^\nu}\left(-\frac{C^2}{4A}\right)e_{p^\nu}\left(-\frac{D^2}{4B}\right),
	\end{aligned}
	\end{equation*}
	with the indicator function
	
	\begin{equation*}
	\chi(p^\nu,A,B,C,D)={\bf 1}\left\{
	\begin{aligned}
        \deg_{p^\nu}(C)\geq \deg_{p^\nu}(A) \\
		\deg_{p^\nu}(D) \geq \deg_{p^\nu}(B)
	\end{aligned}
	\right\}.
	\end{equation*}
 
 If $p=2$, we have
 \begin{equation*}
 S(2^\nu,A,B,C,D)\ll p^{2\nu} \cdot {\bf 1}\{2^{\nu-1}\mid (A,B)\}.
 \end{equation*}
\end{lemma}

The following two lemmas can be found in Zhang \cite{ZHA1}, which will play an important role in the estimation of some integrals.
\begin{lemma}[Non-stationary phase]\label{NSP}
Let $\phi$ be a smooth compactly supported function on $(-\infty,\infty)$ and $f$ be a function which, for $x$ in the support of $\phi$, satisfies:
	\begin{itemize}
		\item $|f'(x)|>A>0,$
		\item $A\geq |f^{(2)}{(x)}|,\cdots,|f^{(N)}(x)|.$
	\end{itemize}
Then
$$\int_{-\infty}^{\infty} \phi(x)e(f(x))dx\ll_{\phi,N}A^{-N}.$$
\end{lemma}

\begin{lemma}[Stationary phase]\label{SP}
	Let $f$ be a quadratic polynomial of two variables $x$ and $y$ whose homogeneous part has discriminant $-D$ with $D>0$. Let $\phi(x,y)$ be a smooth compactly supported function on $\mathbb{R}^2$. Then
	$$\int_{-\infty}^{\infty}\int_{-\infty}^{\infty} \phi(x,y)e(f(x,y))dxdy\ll_{\phi}\frac{1}{\sqrt{D}}.$$
\end{lemma}

The following lemma shows that in $B_T$, matrices can almost be determined by their entries in the bottom row.

\begin{lemma}\label{samebottom}
	Fix $\gamma \in B_T$. We have
	\begin{align*}
	\sum_{\gamma'\in B_T}{\bf 1}\left\{
	\begin{aligned}
	c_1+c_2\sqrt{-\mathfrak D}=c'_1+c'_2\sqrt{-\mathfrak D} \\
	d_1+d_2\sqrt{-\mathfrak D}=d'_1+d'_2\sqrt{-\mathfrak D}
	\end{aligned}
	\right\}\ll 1.
	\end{align*}
\end{lemma} 
\begin{proof}
	Fix $c_1+c_2\sqrt{-\mathfrak D}$ and $d_1+d_2\sqrt{-\mathfrak D}$. The solution of the equation
	$$(a_1+a_2\sqrt{-\mathfrak D})(d_1+d_2\sqrt{-\mathfrak D})-(b_1+b_2\sqrt{-\mathfrak D})(c_1+c_2\sqrt{-\mathfrak D})=1,$$
	can be written as the form
	\begin{align*}
	a_1+a_2\sqrt{-\mathfrak D}&=m_1+n_2\sqrt{-\mathfrak D}+t(c_1+c_2\sqrt{-\mathfrak D}),\\
	b_1+b_2\sqrt{-\mathfrak D}&=m'_1+n'_2\sqrt{-\mathfrak D}+t(d_1+d_2\sqrt{-\mathfrak D}),
	\end{align*}
	where $m_1+n_2\sqrt{-\mathfrak D}$ and $m'_1+n'_2\sqrt{-\mathfrak D}$ is one of these solutions and $t\in \mathbb{Z}[\sqrt{-\mathfrak D}]$.
	
	Since $\parallel a_1+a_2\sqrt{-\mathfrak D} \parallel \leq T$ and $\parallel c_1+c_2\sqrt{-\mathfrak D} \parallel \geq T/10$, there are at most $O(1)$ choices for $t$.
\end{proof}

The following lemma counts the number of solutions of some congruence equations.
\begin{lemma}\label{solcong}
	Let $F,J,K$ be integers, $x\leq A$, $y\leq B$ and $\ell \leq \min(A, B)$. Let $\tau=(F,J)$.
	The number of solutions of $(x,y)$ for 
	$$Fx+Jy\equiv K\quad (\mod \ell)$$
	is 
	\begin{equation*}
	\ll \frac{AB}{\ell/(\ell,\tau)}+A+B.
	\end{equation*}
\end{lemma}
\begin{proof}
	Let $F=F'\tau$ and $J=J'\tau$. If $(\ell, \tau)\mid K$, the congruence equation is equivalent to
	\begin{equation}\label{FJred}
	F'x+J'y\equiv \frac{K}{(\ell,\tau)}\overline{\tau/(\ell,\tau)} \quad (\mod \frac{\ell}{(\ell,\tau)}),
	\end{equation}
	where $\overline{\tau/(\ell,\tau)}$ is the inverse of $\tau/(\ell,\tau)$ modulo $\ell/(\ell,\tau)$.
	Then $-J'y+\frac{K}{(\ell,\tau)}\overline{\tau/(\ell,\tau)}$ is a multiple of $(F',\ell/(\ell,\tau))$. That is
	$$J'y\equiv \frac{K}{(\ell,\tau)}\overline{\tau/(\ell,\tau)} \quad (\mod (F',\frac{\ell}{(\ell,\tau)})).$$
	Thus, since $(F',J')=1$, there are at most
	$$\frac{B}{(F',\frac{\ell}{(\ell,\tau)})}+1$$
	many choices for $y$.
	Fix such a $y$, the equation \eqref{FJred} can be reduced to
	$$F'x\equiv K' \quad (\mod \ell/(\ell,\tau)).$$
	The number of solutions of $x$ is 
	$$\ll \frac{A}{\frac{\ell}{(\ell,\tau)}/(\frac{\ell}{(\ell,\tau)},F')}+1.$$
\end{proof}	

\begin{lemma}\label{countarch}
	Fix $N/2\leq n\leq N$, $1<K\leq T^{1/10}$, and $x,y\asymp X$. Then
	\begin{align*}
	\sum_{\gamma\in B_T}{\bf 1}\left\{|Q_\gamma(x,y)-n|<\frac{N}{K}\right\}\gg \frac{T^{2\delta}}{K}+T^{2\Theta_2},
	\end{align*}
	where $\Theta_2<\delta$ depends only on the spectral gap for $\Gamma$. The implied constant is independent of $x,y$ and $n$.
\end{lemma}
\begin{proof}
	The proof is similar as the argument of Lemma 3.5 in \cite{KZ23}.
\end{proof}

\section{Major Arc Analysis \label{major}}
This section is devoted to deal with the major arcs defined in \eqref{Major} in the circle method. The main conclusion is Proposition \ref{main arcs}. 

From \eqref{T} and \eqref{Major}, the main term is
\begin{equation}\label{MJ}
\begin{aligned}
\mathcal{M}_N(n)=&\int_0^1\sum_{q\leq Q_0}{\sum_{r (\mod q)}}'\sum_{m\in \mathbb{Z}}\mathfrak{t}\left(\frac{N}{K_0}(\theta+m-\frac{r}{q})\right)\widehat{R}(\theta)e(-n\theta)d\theta\\
=& \sum_{q\leq Q_0}{\sum_{r (\mod q)}}' \int_{\mathbb{R}}\mathfrak{t}\left(\frac{N}{K_0}\beta \right)\widehat{R}(\frac{r}{q}+\beta)e(-n(\frac{r}{q}+\beta))d\beta\\
=&\sum_{x,y\in \mathbb{Z}}\Psi(\frac{x}{X})\Psi(\frac{y}{X})\sum_{q\leq Q_0}{\sum_{r (\mod q)}}' \sum_{\gamma\in B_T}\int_{\mathbb{R}}\mathfrak{t}\left(\frac{N}{K_0}\beta \right)e\left((Q_\gamma(x,y)-n)(\frac{r}{q}+\beta)\right)d\beta.
\end{aligned}
\end{equation}

Now we split the set $\gamma\in B_T$ into left cosets of $\Gamma(q)$ 
and apply Lemma \ref{lemmaBKS} with $K=K_0$ ($<T^{1/10}$) to obtain
\begin{equation*}
\begin{aligned}
&\sum_{\gamma\in B_T}e\left((Q_\gamma(x,y)-n)(\frac{r}{q}+\beta)\right)=\sum_{\gamma_0\in \Gamma/\Gamma(q)}e\left((Q_{\gamma_0}(x,y)-n)\frac{r}{q}\right)\sum_{\substack{\gamma\in B_T\\\gamma\equiv \gamma_0(\Gamma(q))}}e\left((Q_\gamma(x,y)-n)\beta\right)\\
=&\sum_{\gamma_0\in \Gamma/\Gamma(q)}e\left((Q_{\gamma_0}(x,y)-n)\frac{r}{q}\right)\left(\frac{1}{[\Gamma:\Gamma(q)]}\left(\sum_{\gamma\in B_T}e\left((Q_\gamma(x,y)-n)\beta\right)\right)+O(T^{2\Theta_1}K_0)\right).
\end{aligned}
\end{equation*}
Plugging the above equation into \eqref{MJ}, we obtain
\begin{align}\label{singular}
\mathcal{M}_N(n)=\sum_{x,y\in \mathbb{Z}}\Psi(\frac{x}{X})\Psi(\frac{y}{X}) \mathfrak{S}_{Q_0; x,y}(n)
\mathfrak M_{x,y}(n)+O\left(\frac{T^{2\Theta_1}X^2K_0^2Q_0^8}{N}\right),
\end{align}
where
\begin{equation*}
\mathfrak{S}_{Q_0;x,y}(n):=\sum_{q\leq Q_0}\frac{1}{[\Gamma:\Gamma(q)]}\sum_{\gamma_0\in \Gamma/\Gamma(q)}{\sum_{r (\mod q)}}' e\left((Q_{\gamma_0}(x,y)-n)\frac{r}{q}\right),
\end{equation*}
and
\begin{equation*}
\mathfrak{M}(n)=\mathfrak{M}_{x,y}(n):=\frac{K_0}{N}\sum_{\gamma\in B_T}\hat{\mathfrak{t}}\left(\frac{K_0}{N}(Q_{\gamma}(x,y)-n)\right).
\end{equation*}

Using Lemma \ref{countarch} and assuming $K_0<\min \{T^{1/10},T^{2(\delta-\Theta_2)}\}$, the Archimedean piece $\mathfrak{M}(n)$ can be bounded from below by
\begin{equation}\label{lowerbdsi}
\mathfrak{M}(n)\gg \frac{T^{2\delta}}{N}.
\end{equation}

To understand $\mathfrak{S}_{Q_0}(n)$, we define a formal singular series
\begin{equation*}
\mathfrak{S}(n)=\sum_{q=1}^\infty U_q(n),
\end{equation*}
with
\begin{equation*}
U_q(n)=\frac{1}{[\Gamma:\Gamma(q)]}\sum_{\gamma_0 \in \Gamma/\Gamma(q)} c_q\left(Q_{\gamma_0}(x,y)-n\right).
\end{equation*}
Here $c_q$ is the Ramanujan sum defined by
$$c_q(m)={\sum_{r(q)}}'e\left(\frac{rm}{q}\right).$$

For fixed $m$, $c_q$ is multiplicative with respect to $q$, and locally
\begin{align}\label{ramanujansumbd}
c_{p^k}(m)=
\begin{cases}
0,& p^\nu\parallel n,\nu\leq k-2,\\
-p^{k-1},& p^{k-1}\parallel m,\\
p^{k-1}(p-1),& p^k\mid m.
\end{cases}
\end{align}

Rewrite $U_q(n)$ as 
\begin{align*}
U_q(n)=\sum_{r(q)}\tau_q(r)c_q(r-n),
\end{align*}
where
\begin{align*}
\tau_q(r)=\frac{1}{[\Gamma:\Gamma(q)]}\sum_{\gamma_0 \in \Gamma/\Gamma(q)}{\bf 1}\{Q_{\gamma_0}(x,y)\equiv r(\mod q)\}.
\end{align*}

As $\Gamma$ is a Zariski dense subgroup of $SL_2(\mathbb Z[\sqrt{-\mathfrak{D}}])$, the strong approximation property for $\Gamma$ implies

\begin{lemma}\label{liftlemma}
	There is an integer $P_{\text{bad}}\geq 1$ such that
	\begin{enumerate}
		\item[i.)] for any $p\nmid P_{\text{bad}}$ and $k\geq 1$,
		\begin{equation}\label{lift}
		\Gamma(p^{k-1})/\Gamma(p^k)=\text{SL}_2(\mathbb Z[\sqrt{-\mathfrak{D}}])(p^{k-1})/\text{SL}_2(\mathbb Z[\sqrt{-\mathfrak{D}}])(p^k);
		\end{equation}
		\item[ii.)] for each of the finitely many primes $p\mid P_{\text{bad}}$, there exists $k_p$ such that \eqref{lift} holds for any $k\geq k_p$.
	\end{enumerate}
\end{lemma}
\begin{proof}
This is the strong approximation of $\text{SL}_2$ (see \cite{Wei}). 
\end{proof}

Using Lemma \ref{liftlemma}, we can prove the following result.
\begin{lemma}
	There is an integer $L\geq 1$ such that
	\begin{enumerate}
		\item[(i)] for any $p\nmid L$ and $k\geq 2$, $U_{p^k}(n)=0$;
		\item[(ii)] for each of the finitely many primes $p\mid L$, there exists $k'_p$ such that $U_{p^k}(n)=0$ for any $k\geq k'_p$.
	\end{enumerate}	
\end{lemma}
\begin{proof} We assume $p$ is odd. We can write
	\begin{equation}\label{gamma0}
	\begin{aligned}
	U_{p^k}(n)=&\frac{1}{[\Gamma:\Gamma(p^k)]}\sum_{\gamma \in \Gamma/\Gamma(p^k)} c_{p^k}\left(Q_{\gamma}(x,y)-n\right)\\
	=&\frac{1}{[\Gamma:\Gamma(p^k)]}\sum_{\gamma_0 \in \Gamma/\Gamma(p^{k-1})}\sum_{\substack{\gamma \in \Gamma(p^{k-1})/\Gamma(p^k)\\\gamma\equiv \gamma_0 (p^{k-1})}} c_{p^k}\left(Q_{\gamma}(x,y)-n\right).
	\end{aligned}
	\end{equation}
	
	For any $\gamma$ with lower row $(c_1+c_2 \sqrt{-\mathfrak{D}}, d_1+d_2 \sqrt{-\mathfrak{D}})$, recalling the definition of $Q_{\gamma}(x,y)$ in \eqref{defQ}, $(\partial_1 Q_\gamma(x,y)/\partial c_1, \partial_1 Q_\gamma(x,y)/\partial c_2, \partial_1 Q_\gamma(x,y)/\partial d_1,\partial_1 Q_\gamma(x,y)/\partial d_2)\neq {\bf{0}}(\mod p)$ (the determinant of the matrix of coefficients of linear equations is $1$).  Together with Lemma \ref{liftlemma}, this leads to that 
	for any $k\geq k_p$, and for any $\gamma_0$ satisfying $Q_{\gamma_0}(x,y)\equiv n \ (\mod p^{k-1})$,
	\begin{equation}\label{problift}
	{\rm Prob}\left\{\gamma\in \Gamma(p^{k-1})/\Gamma(p^k),\gamma\equiv \gamma_0 (p^{k-1})\mid Q_{\gamma}(x,y)\equiv n\ (\mod p^{k})  \right\}=\frac{1}{p}.
	\end{equation}
	According to \eqref{ramanujansumbd}, the innermost sum of \eqref{gamma0} is zero if $Q_{\gamma_0}(x,y)\not\equiv n \ (\mod p^{k-1})$. When $Q_{\gamma_0}(x,y)\equiv n \ (\mod p^{k-1})$ and $k\geq k_p$, the sum is still zero according to \eqref{ramanujansumbd} and \eqref{problift}. Thus, $U_{p^k}(n)=0$ for $k\geq k_p$.
	
\end{proof}

\begin{lemma}\label{bdtau}
	The density function $\tau_p(n)$ has the following property: 
	\begin{align*}
	\tau_{p}(n)=
	\begin{cases}
	\frac{1}{p}+O\left(\frac{1}{p^2}\right),& n\equiv 0\ (\mod p),\\
	\frac{1}{p}+O\left(\frac{1}{p^3}\right),& n\not\equiv 0\ (\mod p),
	\end{cases}
	\end{align*}
	where the implied constants are independent of $n$.
\end{lemma}
\begin{proof}

We only need to consider good odd primes, as the rest primes are finitely many and don't affect the statement.  
	Let $p$ be an odd prime not diving $P_{\rm bad}$, and let 
	\begin{align*}
	V_p=\{(c_1,c_2,&d_1,d_2)\in (\mathbb Z/p\mathbb Z)^4: \\
	&(c_1+c_2\sqrt{-\mathfrak{D}}, d_1+d_2\sqrt{-\mathfrak{D}})\text{ is a row vector of some matrix in } \text{SL}_2(\mathbb Z[\sqrt{-\mathfrak{D}}]/(p))\}.
	\end{align*}
	Let $$W_p=\{(c_1,c_2,d_1,d_2)\in (\mathbb Z/p\mathbb Z)^4: A_{\gamma}x^2+B_{\gamma}y^2+C_{\gamma}x+D_{\gamma}y+E_{\gamma} \equiv n(\mod p)\}.$$
	Then we have 
	\begin{align*}
	\#W_p=&\frac{1}{p}\sum_{a=0}^{p-1}\sum_{c_1,c_2,d_1,d_2\in \mathbb Z/p\mathbb Z}e_p(a(A_{\gamma}x^2+B_{\gamma}y^2+C_{\gamma}x+D_{\gamma}y+E_{\gamma}-n))\\
	=&p^3+\frac{1}{p}\sum_{a=1}^{p-1}\sum_{c_1,c_2,d_1,d_2 \in \mathbb Z/p\mathbb Z}e_p(a(A_{\gamma}x^2+B_{\gamma}y^2+C_{\gamma}x+D_{\gamma}y+E_{\gamma}-n)).
	\end{align*}
	By using Lemma \ref{gauss sum}, we obtain
	\begin{align*}
	\#W_p=p^3+\frac{1}{p}\sum_{a=1}^{p-1}p^2e_p(-an)=	\begin{cases}
	p^3-p,& p\nmid n,\\
	p^3+p(p-1),& p|n.
	\end{cases}
	\end{align*}
Also we have
	\begin{align*}
	\# V_p=
	\begin{cases}
	(p^2-1)^2,& \left(\dfrac{-\mathfrak{D}}{p}\right)=1,\\
	p^4-1,& \left(\dfrac{-\mathfrak{D}}{p}\right)=-1.
	\end{cases}
	\end{align*}

	If $p\nmid n$, any vector in $(c_1,c_2,d_1,d_2)\in W_p$ must have $p\mid (A_\gamma, B_\gamma, C_\gamma, D_\gamma)$, so $(c_1,c_2,c_3,c_4)\in V_p$.  If $p|n$, then 
	\begin{align*}
	\#{(W_p-V_p)}
	=&	\begin{cases} 
	(2p-1)^2,&
	\left(\dfrac{-\mathfrak{D}}{p}\right)=1,\\
	\\
	1,& \left(\dfrac{-\mathfrak{D}}{p}\right)=-1.
	\end{cases}
	\end{align*}

%
%
	For each vector ${\bf{v}}=(c_1,c_2,d_1,d_2)\in V_p$, there are equal numbers of ${\bf{w}}=(a_1,a_2,b_1,b_2)\in (\mathbb Z/p\mathbb Z)^4$ so that $\left(\begin{array}{ccc}  
{\bf w} \\  
{\bf v}
\end{array}
\right)\in \text{SL}(2,\mathbb Z[i]/(p))$.  Therefore, 
	\begin{align*}
	\tau_p(n)= \frac{\#(W_p\cap V_p)}{V_p} 
	=	\begin{cases} \frac{p^3-3p^2+3p-1}{p^4-(2p-1)^2}&
	 \left(\dfrac{-\mathfrak{D}}{p}\right)=1\text{ and } p|n,\\
	\\
	\frac{p+1}{p^2+1},& \left(\dfrac{-\mathfrak{D}}{p}\right)=-1 \text{ and } p|n, \\ \\
		 \frac{p^3-p-1}{p^4-(2p-1)^2} &
              \left(\dfrac{-\mathfrak{D}}{p}\right)=1\text{ and } p\nmid n,\\
	\\
	\frac{p}{p^2+1},& \left(\dfrac{-\mathfrak{D}}{p}\right)=-1 \text{ and } p\nmid n.
	\end{cases}
	\end{align*}
	
	The lemma thus follows.  

\end{proof}

\begin{lemma}\label{bdDp}
	Suppose $n$ is admissible. Let $p$ be an odd prime not dividing $L$. Then $U_p(n)=O\left(1/p^2\right)$ if $p\nmid n$ and $U_p(n)=O\left(1/p\right)$ if $p\mid n$, where the implied constants are independent of $n$.
\end{lemma}
\begin{proof}
	We have
	\begin{align*}
	U_p(n)=&\sum_{r(p)}\tau_p(r)c_p(r-n)= \tau_p(n)(p-1)-\sum_{\substack{r(p)\\ r\not\equiv n(p)}}\tau_p(r)\\
	=&\tau_p(n)(p-1)-(1-\tau_p(n))=p\tau_p(n)-1.
	\end{align*}
	It is easy to see the bound for $U_p(n)$ by using Lemma \ref{bdtau}.
\end{proof}

%
%

From the above two lemmas,  it is straightforward to show that,
\begin{lemma}\label{sinser}
	We have
\begin{align*}
\mathfrak{S}_{Q_0;x,y}(n)\gg \frac{1}{\log n}, \qquad \text{$n$ is admissible,}
\end{align*}
and
\begin{align*}
\mathfrak{S}_{Q_0;x,y}(n)\ll \frac{\log n}{Q_0}, \qquad \text{$n$ is not admissible.}
\end{align*}
\end{lemma}
\begin{proof}
The proof is the same as the proof for Proposition 6.1 in \cite{FSZ}.
\end{proof}

Combining \eqref{singular}, \eqref{lowerbdsi} and Lemma \ref{sinser}, we can make the main conclusion of this section.

\begin{proposition}\label{main arcs}
	 For any admissible $n\in [N,2N]$, we have 
	\begin{align*}
	\mathcal{M}_N(n)\gg \frac{T^{2\delta-2}}{\log n},
	\end{align*}
	assuming 
	$$K_0<\min \{T^{1/10},T^{2(\delta-\Theta_2)}\}\  \mbox{and}\ \  K_0^2Q_0^8<T^{2(\delta-\Theta_1)}.$$
\end{proposition}


\section{Minor Arcs Analysis \label{minor}}
We can not give a satisfactory individual bound for $\mathcal{E}_N(n)$. Instead, we give an $\mathcal{L}^2$ bound for $\mathcal E_N$ which implies that $\mathcal{E}_N(n)$ is small in average.
\begin{proposition}\label{minor arcs}
	There exists a positive $\eta$ such that
	$$\sum_{n\in \mathbb{Z}}|\mathcal{E}(n)|^2\ll T^{4\delta-4}N^{1-\eta},$$
	where the implied constant depends only on $\Gamma$ and $\mathfrak{D}$.  
\end{proposition}

Denote by
\begin{equation}\label{integral}
I=\int_0^1\left|1-\mathfrak{T}(\theta)\right|^2\left|\widehat{R}(\theta)\right|^2d\theta.
\end{equation}
By Plancherel, it is enough to prove
$$I\ll T^{4\delta-4}N^{1-\eta}.$$
Let 
$$M=T^2X.$$ 
Let $K_0$ and $Q_0$ be a small power of $T$ as before.
The integral \eqref{integral} is bounded by the sum of the following four integrals:
\begin{align*}
I_1&=\sum_{q\leq Q_0}{\sum_{r (\mod q)}}'\int_{-K_0/N}^{K_0/N}\left(\frac{N\beta}{K_0}\right)^2\left|\widehat{R}(\frac{r}{q}+\beta)\right|^2d\beta,\\
I_2&=\sum_{q\leq Q_0}{\sum_{r (\mod q)}}'\int_{K_0/N\leq |\beta|\leq 1/qM}\left|\widehat{R}(\frac{r}{q}+\beta)\right|^2d\beta,\\
I_3&=\sum_{Q_0< q\leq X^{1-\omega}}{\sum_{r (\mod q)}}'\int_{|\beta|\leq 1/qM}\left|\widehat{R}(\frac{r}{q}+\beta)\right|^2d\beta,\\
I_4&=\sum_{X^{1-\omega}<q\leq M}{\sum_{r (\mod q)}}'\int_{|\beta|\leq 1/qM}\left|\widehat{R}(\frac{r}{q}+\beta)\right|^2d\beta.
\end{align*}
Here $\omega$ is a sufficiently small constant depending on $\delta$ (see \eqref{II=Ineq} for the explanation).

Recalling \eqref{defRhat}, we have
\begin{align*}
\widehat{R}(\frac{r}{q}+\beta)
&= \sum_{\gamma \in B_T}\sum_{x_0,y_0 (\mod q)}e\left(Q_\gamma(x_0,y_0)\frac{r}{q}\right)\sum_{\substack{x\equiv x_0 (q)\\y\equiv y_0 (q)}}\Psi\left(\frac{x}{X}\right)\Psi\left(\frac{y}{X}\right)e\left(Q_\gamma(x,y)\beta\right).
\end{align*}
According to the definition of \eqref{defQ}, we can rewrite $\widehat{R}(\frac{r}{q}+\beta)$ as the form
\begin{align*}
\widehat{R}(\frac{r}{q}+\beta)=&\sum_{\gamma \in B_T}e\left(E_\gamma (\frac{r}{q}+\beta) \right)\sum_{x_0,y_0 (\mod q)}e\left(\frac{A_\gamma x_0^2+C_\gamma x_0}{q}r\right)e\left(\frac{B_\gamma x_0^2+D_\gamma x_0}{q}r\right)\\
                                              &\times \sum_{x\equiv x_0(q)}\Psi\left(\frac{x}{X}\right)e((A_\gamma x^2+C_\gamma x)\beta)\sum_{y\equiv y_0(q)}\Psi\left(\frac{y}{X}\right)e((B_\gamma y^2+D_\gamma y)\beta).
\end{align*}
Then by Poisson summation,
\begin{align}\label{Rhat}
\widehat{R}(\frac{r}{q}+\beta)=&\sum_{\gamma \in B_T}\mathfrak{R}_\gamma(\frac{r}{q}+\beta)e(E_\gamma (\frac{r}{q}+\beta))
\end{align}
with
\begin{align}\label{defR}
\mathfrak{R}_\gamma(\frac{r}{q}+\beta)=X^2\sum_{\xi,\zeta\in \mathbb{Z}}S_{\gamma}(q,r,\xi,\zeta)J_{\gamma}(\beta,q,\xi,\zeta),
\end{align}
where
\begin{align}\label{defs1}
S_{\gamma}(q,r,\xi,\zeta)=\frac{1}{q^2}\sum_{x_0,y_0(q)}e\left(\frac{A_\gamma x_0^2+C_\gamma x_0}{q}r-\frac{x_0}{q}\xi\right)e\left(\frac{B_\gamma y_0^2+D_\gamma y_0}{q}r-\frac{y_0}{q}\zeta\right),
\end{align}
and
\begin{align*}
J_{\gamma}(\beta,q,\xi,\zeta)=\int_{\mathbb{R}}\int_{\mathbb{R}}&\Psi(x)\Psi(y)e\left(A_\gamma \beta X^2x^2+C_\gamma\beta Xx+\frac{\xi X}{q}x\right)\\
&\times e\left(B_\gamma \beta X^2y^2+D_\gamma\beta Xy+\frac{\zeta X}{q}y\right)dxdy.
\end{align*}

Lemma \ref{gauss sum}, for $q=p^\nu$,
implies the bound for $S_\gamma(q,r,\xi,\zeta)$,

\begin{align}\label{trivialbdS}
S_\gamma(q,r,\xi,\zeta)\ll
\begin{cases}
0,& q\mid A_\gamma, B_\gamma, q\nmid (rC_\gamma-\xi,rD_\gamma-\zeta),\\
1,& q\mid A_\gamma, B_\gamma, q\mid (rC_\gamma-\xi,rD_\gamma-\zeta),\\
q^{-1}(q,A_\gamma)^{1/2}(q,B_\gamma)^{1/2},& q\nmid A_\gamma, \ (q\mid B_\gamma \ \text{or}\  q\nmid B_{\gamma}).
\end{cases}
\end{align}

Trivially we have the bound for $J_\gamma(\beta,q,\xi,\zeta)$,
$$J_\gamma(\beta,q,\xi,\zeta)\ll 1.$$
For $|\beta|<\frac{1}{qM}$ and $\xi\in\mathbb{Z}$, we obtain
\begin{align*}
A_\gamma\beta X^2x\ll& \frac{\xi X}{q}, \quad \text{if $\xi\neq 0$, $0.5<x<2.5$,}\\
C_\gamma \beta X \ll& \frac{\xi X}{q}, \quad \text{if $\xi\neq 0$}.
\end{align*}
Thus, non-stationary phase Lemma \ref{NSP} implies, for any absolute $L_1,L_2>0$,
\begin{align}\label{nonstationarybd}
J_{\gamma}(\beta,q,\xi,\zeta)\ll \left(\frac{\xi X}{q}\right)^{-L_1}\left(\frac{\zeta X}{q}\right)^{-L_2}, \quad \text{if $\zeta,\xi\neq 0$, $q\ll X$, $|\beta|<\frac{1}{qM}$}.
\end{align}
Also, stationary phase Lemma \ref{SP} implies
\begin{align}\label{stationarybd}
J_\gamma(\beta,q,\xi,\zeta)\ll \frac{1}{\sqrt{A_\gamma \beta X^2}}\frac{1}{\sqrt{B_\gamma \beta X^2}}=\frac{1}{A_\gamma\mathfrak{D}^{1/2} \beta X^2} \ll \frac{1}{\beta T^2X^2}.
\end{align}

For $\gamma,\gamma'\in B_T$, denote
\begin{align*}
S(q,\gamma,\xi,\zeta,\gamma',\xi',\zeta')={\sum_{r(q)}}'S_\gamma(q,r,\xi,\zeta)\overline{S_{\gamma'}(q,r,\xi',\zeta')}e(E_\gamma\frac{r}{q})e(-E_{\gamma'}\frac{r}{q}).
\end{align*}
Then $S(q,\gamma,\xi,\zeta,\gamma',\xi',\zeta')$ is multiplicative over $q$.

Let $p$ be an odd prime ($p\mid \mathfrak{D}$ or not) and $p^\nu \nmid A_\gamma, A'_\gamma$. Let $\ell_0=(B_\gamma,p^\nu)$ and $\ell'_0=(B_{\gamma'},p^\nu)$. By Lemma \ref{ESQF},
%
we can write (the process is a little complicated by considering the cases $p\mid \mathfrak{D}$ or not)
\begin{align}\label{s1s2}
S(p^\nu,\gamma,\xi,\zeta,\gamma',\xi',\zeta')=S_1\cdot S_2,
\end{align}
with
\begin{equation}\label{s1}
	\begin{aligned}
		S_1=&\frac{1}{p^{2\nu}}(p^\nu,A_\gamma)^{1/2}(p^\nu,B_\gamma)^{1/2}(p^\nu,A_{\gamma'})^{1/2}(p^\nu,B_{\gamma'})^{1/2}e\left(\mathfrak{A}\right)\\
		&\times e\left(-\mathfrak{A}'\right)i^{\epsilon\left(\frac{p^\nu}{(p^\nu,A_\gamma)}\right)+\epsilon\left(\frac{p^\nu}{(p^\nu,A_\gamma')}\right)+\epsilon\left(\frac{p^\nu}{(p^\nu,B_\gamma)}\right)+\epsilon\left(\frac{p^\nu}{(p^\nu,B_{\gamma'})}\right)},
	\end{aligned}
\end{equation}
\begin{equation}\label{s2}
\begin{aligned}
S_2
      =&{\sum_{r (p^\nu)}}' e_{p^\nu}(A_\gamma r)e_{p^\nu}(-A_{\gamma'}r)e\left(\frac{-\mathfrak{B}+\mathfrak{B}'}{p^\nu}\right)\chi(r),
\end{aligned}
\end{equation}
where
$$\mathfrak{A}=\frac{\overline{2}\overline{B_\gamma /\ell_0}(D_\gamma \zeta+\mathfrak{D}C_\gamma \xi)/\ell_0}{p^\nu},\quad  \mathfrak{A}'=\frac{\overline{2}\overline{B_{\gamma'} /\ell'_0}(D_{\gamma'} \zeta'+\mathfrak{D}C_{\gamma'} \xi')/\ell'_0}{p^\nu},$$ 
and
$$\mathfrak{B}=\overline{4r}\overline{B_\gamma/\ell_0}\frac{(\zeta^2+\mathfrak{D}\xi^2)/\ell_0}{p^\nu}
,\quad  \mathfrak{B}'=\overline{4r}\overline{B_{\gamma'}/\ell'_0}\frac{(\zeta'^2+\mathfrak{D}\xi'^2)/\ell'_0}{p^\nu},$$
with $\overline{B_\gamma /\ell_0}$ (similarly for $\overline{B_{\gamma'} /\ell'_0}$) being the one satisfying $\overline{B_\gamma /\ell_0} \cdot B_\gamma/\ell_0 \equiv 1 (\mod p^\nu)$ (since $p^\nu \nmid A_\gamma$, $\gcd(B_\gamma/\ell_0, p^\nu)=1$).
Here
$\chi(\cdot)$ is the indicator function defined by
\begin{equation}\label{defchi}
\chi(\cdot)={\bf 1}\left\{
\begin{aligned}
\deg_{p^\nu} (rC_\gamma-\xi)\geq \deg_{p^\nu}(A_\gamma)\\
\deg_{p^\nu} (rD_\gamma-\zeta)\geq \deg_{p^\nu}(B_\gamma)
\end{aligned}
\right\}\cdot {\bf 1}\left\{
\begin{aligned}
\deg_{p^\nu} (rC_{\gamma'}-\xi)\geq \deg_{p^\nu}(A_{\gamma'})\\
\deg_{p^\nu} (rD_{\gamma'}-\zeta)\geq \deg_{p^\nu}(B_{\gamma'})
\end{aligned}
\right\}.
\end{equation}

\begin{remark}
$\bullet$ Here we have used the fact $\mathfrak{D}A_\gamma=B_{\gamma}$ and $$\mathfrak{D}C_\gamma^2+D_\gamma^2=4\mathfrak{D}(c_1^2+c_2^2\mathfrak{D})(d_1^2+d_2^2\mathfrak{D})=4\mathfrak{D}(d_1^2+d_2^2\mathfrak{D})A_\gamma=4\mathfrak{D}(E_{\gamma}-A_\gamma)A_\gamma.$$

$\bullet$ It is easy to  see $\chi(*)$ is an indicator function of an arithmetic progression set and when $\chi(\cdot)=1$
$$\deg_{p^\nu}(\mathfrak{D}\xi^2+\zeta^2)\geq \deg_{p^\nu}B_\gamma.$$
\end{remark}


\begin{lemma} Let $p$ be an odd prime and $p^\nu \nmid A_\gamma, A_{\gamma'}$.
When $A_{\gamma}\neq A_{\gamma'}$, we have
\begin{align*}
S(p^\nu,\gamma,\xi,\zeta,\gamma',\xi',\zeta')\ll  \frac{(p^\nu,A_\gamma)(p^\nu,A_{\gamma'})}{p^{2\nu}} p^{3\nu/4+\ve}(p^\nu, A_{\gamma}-A_{\gamma'})^{1/4}.
\end{align*}
When $A_{\gamma}= A_{\gamma'}$ and $\ell_0=(p^\nu,\mathfrak{D}A_\gamma)$, we have
\begin{align*}
S(p^\nu,\gamma,\xi,\zeta,\gamma',\xi',\zeta')\ll \frac{(p^\nu,A_\gamma)(p^\nu,A_{\gamma'})}{p^{2\nu}} p^{3\nu/4+\ve}\left(p^\nu, {\overline{4\mathfrak{D}A_\gamma/\ell_0}(\mathfrak{D}\xi^2+\zeta^2-\mathfrak{D}\xi'^2-\zeta'^2)/\ell_0}\right)^{1/4}.
\end{align*}
\end{lemma}
\begin{proof}
	If $p^\nu \nmid A_{\gamma}$ and $p^\nu \nmid A_{\gamma'}$, applying Kloosterman's elementary $3/4$ bound for this type of sum (see Lemma 3.9 in \cite{ZHA1}), we can get the bound easily. 
	\end{proof}

 Using the multiplicity of $S(q,\gamma,\xi,\zeta,\gamma',\xi',\zeta')$ over $q$, we can obtain the following bound.

\begin{corollary}\label{corSodd}
Let $q=p_1^{\nu_1}\cdots p_s^{\nu_s}$ be an odd integer with $p_i^{\nu_i}\nmid A_\gamma, A_{\gamma'}$ for $i=1,2,\cdots, s$. When $A_{\gamma}\neq A_{\gamma'}$, we have
\begin{align}\label{bdSneq}
S(q,\gamma,\xi,\zeta,\gamma',\xi',\zeta')\ll (q,A_{\gamma})(q,A_{\gamma'})(q, A_{\gamma}-A_{\gamma'})^{1/4}q^{-5/4+\ve}.
\end{align}

When $A_{\gamma}= A_{\gamma'}$ and $\ell_0=(q,\mathfrak{D}A_\gamma)$, we have
\begin{align}\label{=neq}
S(q,\gamma,\xi,\zeta,\gamma',\xi',\zeta')\ll (q,A_{\gamma})^2\left(q, {\overline{4\mathfrak{D}A_\gamma/\ell_0}(\mathfrak{D}\xi^2+\zeta^2-\mathfrak{D}\xi'^2-\zeta'^2)/\ell_0}\right)^{1/4}q^{-5/4+\ve}.
\end{align}

\end{corollary}

Besides, for $q=2^\nu$ or $q\mid {\rm lcm}(A_\gamma, A_{\gamma'})$, by using Lemma \ref{gauss sumeven} and \eqref{trivialbdS}, we have
\begin{corollary}\label{corSeven}
	For $q_1=2^\nu$ and $q_2\mid {\rm lcm}(A_\gamma, A_{\gamma'})$, we have 
	\begin{equation}\label{bd2power0} 
	S(2^\nu,\gamma,\xi,\zeta,\gamma',\xi',\zeta') \ll 2^\nu \cdot {\bf 1}\{2^{\nu-1}\mid (A_\gamma,A_{\gamma'})\},
	\end{equation}
and
\begin{equation}\label{bd2power1} 
	S(q_2,\gamma,\xi,\zeta,\gamma',\xi',\zeta') \ll \frac{(q_2,A_\gamma)(q_2,A_{\gamma'})}{q_2}.
\end{equation}
\end{corollary}

\section{Estimate for $I_1$ and $I_2$}
In this section, we shall give upper bounds for $I_1$ and $I_2$.
$$I_1=\sum_{q\leq Q_0}{\sum_{r(q)}}'\int_{-K_0/N}^{K_0/N}\left(\frac{N\beta}{K_0}\right)^2|\widehat{R}(\frac{r}{q}+\beta)|^2d\beta.$$


According to \eqref{defR}, \eqref{trivialbdS}, \eqref{nonstationarybd} and \eqref{stationarybd}, main contribution comes from $\xi,\zeta=0$. We have
\begin{align*}
\mathfrak{R}_\gamma(\frac{r}{q}+\beta)\ll X^2 \frac{1}{A_\gamma \beta X^2}\ll \frac{1}{T^2\beta}.
\end{align*}
Then, by \eqref{Rhat},
\begin{align}\label{bdR0}
\widehat{R}(\frac{r}{q}+\beta)\ll T^{2\delta-2}\frac{1}{\beta},
\end{align}
which implies
\begin{align*}
I_1\ll \sum_{q\leq Q_0} q T^{4\delta-4}\frac{N^2}{K_0^2}\frac{K_0}{N} \ll T^{4\delta-4}N\frac{Q_0^2}{K_0}.
\end{align*}
Then 
\begin{equation}\label{bdI1}
I_1\ll T^{4\delta-4}N^{1-\eta_0},
\end{equation}
if $K_0\gg Q_0^2.$

The estimate for $I_2$ is similar.

$$I_2=\sum_{q\leq Q_0}{\sum_{r(q)}}'\int_{K_0/N}^{1/qM}|\widehat{R}(\frac{r}{q}+\beta)|^2d\beta.$$
By using \eqref{bdR0} again, we have
\begin{align*}
I_2\ll \sum_{q\leq Q_0} q T^{4\delta-4}\frac{N}{K_0} \ll T^{4\delta-4}N \frac{Q^2_0}{K_0}.
\end{align*}
Then 
\begin{equation}\label{bdI2}
I_2\ll T^{4\delta-4}N^{1-\eta_0},
\end{equation}
if $K_0\gg Q_0^2.$


\section{Estimate for $I_3$}
For $Q_0<Q \leq X^{1-\omega}$, denote by
\begin{align}\label{defIQ}
I_Q=\sum_{Q< q\leq 2Q}{\sum_{r (\mod q)}}'\int_{-1/qM}^{1/qM}\left|\widehat{R}(\frac{r}{q}+\beta)\right|^2d\beta.
\end{align}

Observing that
\begin{align*}
{\sum_{r (\mod q)}}'\left|\widehat{R}(\frac{r}{q}+\beta)\right|^2 =&{\sum_{r (\mod q)}}' X^4 \sum_{\xi,\zeta,\xi',\zeta'\in \mathbb{Z}}\sum_{\gamma,\gamma'}e(E_\gamma\theta-E_{\gamma'}\theta)J_{\gamma}(\beta,q,\xi,\zeta)\overline{J_{\gamma'}(\beta,q,\xi',\zeta')}\\
&\times S_\gamma(q,r,\xi,\zeta)\overline{S_{\gamma'}(q,r,\xi',\zeta')}\\
=:&S^{=}+ S^{\neq},
\end{align*}
where
\begin{align*}
S^{=}=X^4 \sum_{\xi,\zeta,\xi',\zeta'\in \mathbb{Z}}\sum_{\substack{\gamma,\gamma'\\A_\gamma=A_{\gamma'}}}e(E_\gamma\beta-E_{\gamma'}\beta)J_{\gamma}(\beta,q,\xi,\zeta)\overline{J_{\gamma'}(\beta,q,\xi',\zeta')} S(q,\gamma,\xi,\zeta,\gamma',\xi',\zeta'),
\end{align*}
\begin{align*}
S^{\neq}=X^4 \sum_{\xi,\zeta,\xi',\zeta'\in \mathbb{Z}}\sum_{\substack{\gamma,\gamma'\\A_\gamma\neq A_{\gamma'}}}e(E_\gamma\beta-E_{\gamma'}\beta)J_{\gamma}(\beta,q,\xi,\zeta)\overline{J_{\gamma'}(\beta,q,\xi',\zeta')} S(q,\gamma,\xi,\zeta,\gamma',\xi',\zeta'),
\end{align*}
we can write
\begin{align}\label{II=Ineq}
I_Q=I_Q^=+I_Q^{\neq}:=\int_{-1/qM}^{1/qM} \sum_{Q<q\leq 2Q} S^= d\beta+\int_{-1/qM}^{1/qM} \sum_{Q<q\leq 2Q} S^{\neq} d\beta.
\end{align}
 The main contribution comes from $\xi,\zeta,\xi',\zeta'=0$. For other cases, we can have sufficient saving from $J_\gamma$ and $J_{\gamma'}$ by \eqref{nonstationarybd}.  (That's why we choose $q\leq X^{1-\omega}$. We can get around the trouble from terms $q\sim X$.)
 
We estimate $S^=$ first.
Using \eqref{stationarybd} for $J_\gamma$ and the fact
$$\int_{-1/qM}^{1/qM} \min \{1, \frac{1}{T^4X^4\beta^2}\}d\beta \ll \frac{1}{T^2X^2},$$
we can get by \eqref{=neq}-\eqref{bd2power1} that (here $(2,q'q'')=1$, $(q',q'')=1$, $q'=p_1^{\nu_1}\cdots p_s^{\nu_s}$ with $p_i\nmid A_\gamma, A_{\gamma'}$ and $q''\mid {\rm lcm}\{A_\gamma, A_{\gamma'}\}$)
\begin{equation}\label{bdS=}
\begin{aligned}
I_Q^= \ll& X^4 \sum_{\substack{\gamma,\gamma'\\A_\gamma=A_{\gamma'}}}\sum_{\substack{Q<q\leq 2Q\\q=2^\nu q'q''}} \frac{(q',A_\gamma)^2}{q'}q''\frac{1}{T^2X^2} \times 2^\nu\cdot {\bf 1}\{2^{\nu-1}\mid A_\gamma \}\\
\ll& \frac{X^2}{QT^2}\sum_{\substack{\gamma,\gamma'\\A_\gamma=A_{\gamma'}}} \sum_{\substack{2^\nu\leq 2Q\\2^{\nu-1}\mid A_\gamma}}4^\nu \sum_{d\leq T^2, dq''|A_{\gamma}} \sum_{Q<2^\nu q'q''\leq 2Q, d|q'}d^2 q''^2\\
\ll&
\frac{X^2}{T^2} \sum_{\substack{d'\leq T^2\\(d',2)=1}}d' \sum_{\gamma, d'|A_\gamma} \sum_{\gamma', A_{\gamma'}=A_{\gamma}}1\sum_{\substack{2^\nu\leq 2Q\\2^{\nu-1}\mid A_\gamma}}2^\nu,
\end{aligned}
\end{equation}
where we write $d'=dq''$ in the last display.

Write $d'=\ell_1\ell_2$ with $\ell_1\mid \mathfrak{D}$ and $(\ell_2,\mathfrak{D}/\ell_1)=1$. Then
\begin{equation*}
\begin{aligned}
I_Q^= \ll
\frac{X^2}{T^2}\sum_{\ell_1\mid \mathfrak{D}} \sum_{\substack{\ell_1\ell_2\leq T^2\\(\ell_1\ell_2,2)=1}}\ell_1\ell_2 \sum_{\gamma, \ell_1\ell_2|A_\gamma} \sum_{\gamma', A_{\gamma'}=A_{\gamma}}1\sum_{\substack{2^\nu\leq 2Q\\2^{\nu-1}\mid A_\gamma}}2^\nu.
\end{aligned}
\end{equation*}

Since $\xi,\zeta,\xi',\zeta'=0$, by \eqref{defchi}, if $d|(A_{\gamma},q)$, we have $d|C_\gamma, D_\gamma,C_{\gamma'}, D_{\gamma'}$. Consider the problem in $\mathbb{Z}[\sqrt{\mathfrak{-D}}]$. According to the definition of $C_\gamma, D_\gamma$ in \eqref{defQ}, we have
\begin{align*}
\ell_2\mid 2(c_1+c_2\sqrt{\mathfrak{-D}})(d_1-d_2\sqrt{\mathfrak{-D}}).
\end{align*}
Also we have
\begin{align*}
\ell_2\mid (c_1+c_2\sqrt{\mathfrak{-D}})(c_1-c_2\sqrt{\mathfrak{-D}}).
\end{align*}
Since in $\mathbb{Z}[\sqrt{\mathfrak{-D}}]$, $c_1+c_2\sqrt{\mathfrak{-D}}$ and $d_1+d_2\sqrt{\mathfrak{-D}}$ are relatively prime, we can get
\begin{align*}
\ell_2\mid 2(c_1-c_2\sqrt{\mathfrak{-D}}),
\end{align*}
which means
$$\ell_2\mid c_1,c_2.$$

The contribution from terms $\ell_1 \mid \mathfrak{D}$ is negligible up to a constant. We assume $\ell_1=1$ for simplicity. Thus, thanks to Lemma \ref{samebottom}, \eqref{bdS=} can be bounded by
\begin{equation}\label{I=}
\begin{aligned}
I_Q^=\ll& \frac{X^2}{T^2} \sum_{\substack{2^\nu\leq 2Q}}2^\nu\sum_{d\leq T^2} d \sum_{\substack{c_1,c_2\leq T,d|(2c_1,2c_2)\\2^{\nu-1}\mid c_1^2+c_2^2\mathfrak{D}}}\sum_{\substack{d_1,d_2\leq T}}
 \sum_{\substack{d'_1,d'_2 \leq T }}
 \sum_{\substack{c'_1,c'_2\leq T\\{c'}_1^2+{c'}_2^2\mathfrak{D}=c_1^2+c_2^2\mathfrak{D}}}1\\
 \ll& \frac{X^2}{T^2}N^\ve \sum_{d\leq T^2}d\frac{T^2}{d^2}T^4\ll X^2T^4N^\ve\ll T^{4\delta-4}N^{1-\eta_0},
\end{aligned}
\end{equation}
if $\delta>3/2$.


For $I_Q^{\neq}$, using \eqref{bdSneq}-\eqref{bd2power1}, we have the bound (take $\ell_1=(q'',A_{\gamma})$, $\ell_2=(q'',A_{\gamma'})$, $k_1=(q',A_\gamma)$, $k_2=(q',A_{\gamma'})$ and $k_3=(q',A_{\gamma}-A_{\gamma'})$)
\begin{align*}
I_Q^{\neq}
\ll& X^4\sum_{\substack{\gamma,\gamma'\\A_\gamma\neq A_{\gamma'}}}\sum_{\substack{Q<q\leq 2Q\\q=2^\nu q'q''}}q''^{-1}\ell_1\ell_2 \cdot 2^\nu \cdot {\bf 1}\{2^{\nu-1}\mid (A_\gamma,A_\gamma') \}\int_{-1/qM}^{1/qM}\min\{1,\frac{1}{T^4X^4\beta^2}\}d\beta  \\
&\times (q',A_{\gamma})(q',A_{\gamma'})(q',A_{\gamma}-A_{\gamma'})^{1/4}q'^{-5/4+\ve}\\
\ll& \frac{X^2}{T^2Q^{5/4-\ve}}\sum_{2^\nu q'' \leq 2Q} {q''}^{1/4-\varepsilon}(2^\nu)^{9/4-\varepsilon} \sum_{Q<2^\nu q'q''\leq 2Q}\sum_{\substack{\gamma,\gamma',2^{\nu-1}\mid (A_\gamma,A_{\gamma'})\\A_\gamma\neq A_{\gamma'}}}\ell_1\ell_2k_1k_2k_3^{1/4}.
\end{align*}
Denote $h=\gcd(k_1,k_2)$. Then
\begin{align*}
I_Q^{\neq} \ll& \frac{X^2N^\ve}{T^2Q^{5/4}}\sum_{2^\nu q'' \leq 2Q} {q''}^{1/4-\varepsilon}(2^\nu)^{9/4-\varepsilon} \sum_{h\leq Q/2^\nu}\sum_{h|k_1,k_1\leq Q/2^\nu}\sum_{h|k_2,k_2\leq Q/2^\nu}\sum_{\gamma, 2^\nu k_1|2A_\gamma}
\\
&\times \sum_{\gamma', 2^\nu k_2|2A_{\gamma'}}\ell_1\ell_2 k_1k_2 \sum_{k_3\leq Q/2^\nu,h|k_3|A_\gamma- A_{\gamma'}}k_3^{1/4}
\sum_{Q<2^\nu q'q''\leq 2Q}1\{[k_1,k_2,k_3]|q'\}.
\end{align*}
Using the fact that $k_1/h,k_2/h,k_3/h$ are mutually relatively prime and
$$\sum_{Q/q''2^\nu<q\leq 2Q/q''2^\nu}1\{[k_1,k_2,k_3]|q'\} \leq \frac{Qh^2}{q''2^\nu k_1k_2k_3},$$
we can write
\begin{align*}
I_Q^{\neq} 
\ll&\frac{X^2N^\ve}{T^2Q^{5/4}} \sum_{2^\nu q'' \leq 2Q} {q''}^{-3/4-\varepsilon}(2^\nu)^{5/4-\varepsilon}\sum_{h\leq Q/2^\nu}\sum_{h|k_1,k_1\leq Q/2^\nu}\sum_{h|k_2,k_2\leq Q/2^\nu}\\
&\times \sum_{\gamma, 2^\nu q'' k_1|2A_\gamma}\sum_{\gamma', 2^\nu k_2|2A_{\gamma'}}\ell_1\ell_2 k_1k_2\sum_{k_3\leq Q/2^\nu,h|k_3|A_\gamma- A_{\gamma'}}k_3^{1/4}\frac{Qh^2}{k_1k_2k_3}.
\end{align*}
Using similar counting technique as in \eqref{bdS=} and \eqref{I=} and the fact that $\ell_1\ell_2>q''$, we have
\begin{equation}\label{Ineq}
\begin{aligned}
I_Q^{\neq}
\ll& 
\frac{X^2N^\ve}{T^2Q^{1/4}}\sum_{2^\nu q''\leq 2Q} {q''}^{-\tfrac{3}{4}}2^{\tfrac{9\nu}{4}}\sum_{h\ll T^2}h^2\sum_{k_1,k_2\leq T^2/h} \sum_{\gamma,2^\nu hk_1|2A_{\gamma}}\sum_{\substack{\gamma',2^\nu hk_2|2A_{\gamma'}}}\ell_1\ell_2 \sum_{k_3\mid (A_\gamma-A_{\gamma'})}\frac{1}{k_3^{3/4}}\\
\ll& \frac{X^2N^\ve}{T^2Q^{1/4}}\sum_{h\leq Q}h^2\sum_{k_1,k_2\leq T^2/h} \frac{T^4}{(hk_1)^2}\frac{T^4}{(hk_2)^2}
\ll \frac{NT^4N^\ve}{Q^{1/4}}
\ll N^{1-\eta_0}T^{4\delta-4},
\end{aligned}
\end{equation}
when $T^{8-4\delta}\ll Q_0^{1/4}$.

By \eqref{II=Ineq}, \eqref{I=} and \eqref{Ineq}, we can conclude that, if $\delta>3/2$ and $T^{8-4\delta}\ll Q_0^{1/4}$,
\begin{equation}\label{bdI3}
I_3\ll T^{4\delta-4}N^{1-\eta_0}.
\end{equation}


\section{Estimate for $I_4$}
For $X^{1-\omega}<Q\leq M$, we use the same symbol as \eqref{defIQ},
\begin{align*}
I_Q=\sum_{Q< q\leq 2Q}{\sum_{r (\mod q)}}'\int_{-1/qM}^{1/qM}\left|\widehat{R}(\frac{r}{q}+\beta)\right|^2d\beta.
\end{align*}
Recall that
\begin{align*}
\widehat{R}\left(\frac{r}{q}+\beta \right)=&\sum_{\gamma}\sum_{x,y}\Psi(\frac{x}{X})\Psi(\frac{y}{X})e((A_\gamma x^2+C_\gamma x)\beta)e((B_\gamma y^2+D_\gamma y)\beta)e(E_\gamma (\frac{r}{q}+\beta)  )\\
&\times e((A_\gamma x^2+C_\gamma x)\frac{r}{q})e((B_\gamma y^2+D_\gamma y)\frac{r}{q}).
\end{align*}
Inserting extra harmonics by writing $e((A_\gamma x^2+C_\gamma x)\frac{r}{q})$ into its Fourier expansion, 
\begin{equation*}
e((A_\gamma x^2+C_\gamma x)\frac{r}{q})=\frac{1}{q}\sum_{m,s(q)}e((A_\gamma s^2+C_\gamma s)\frac{r}{q}+\frac{sm}{q})e(-\frac{mx}{q}),
\end{equation*}
\begin{equation*}
e((B_\gamma y^2+D_\gamma y)\frac{r}{q})=\frac{1}{q}\sum_{n,t(q)}e((B_\gamma t^2+D_\gamma t)\frac{r}{q}+\frac{tn}{q})e(-\frac{ny}{q}),
\end{equation*}
we can write
\begin{align*}
\widehat{R}\left(\frac{r}{q}+\beta \right)=\sum_{\gamma}\sum_{m,n(q)}S_{\gamma}(q,r,-m,-n)\lambda_{\gamma}(\beta,X,\frac{m}{q},\frac{n}{q})e(E_\gamma\theta).
\end{align*}
Here
$$\lambda_\gamma(\beta,X,s,t)=\sum_{x,y\in \mathbb{Z}}\Psi(\frac{x}{X})\Psi(\frac{y}{X})e((A_\gamma x^2+C_\gamma x)\beta)e((B_\gamma y^2+D_\gamma y)\beta)e(-sx)e(-ty)$$
and
$S_{\gamma}(q,r,m,n)$ is defined by \eqref{defs1}.
Trivially, 
\begin{equation}\label{triviallambda}
\lambda_\gamma(\beta,X,s,t)\ll X^2.
\end{equation}
We write
\begin{align*}
I_Q
=&\sum_{Q<q\leq 2Q}\sum_{\gamma,\gamma'}\sum_{m,n(q)}\sum_{m',n'(q)}S(q,\gamma,-m,-n,\gamma',-m',-n')\\
&\times \int_{-1/qM}^{1/qM}\lambda_{\gamma}(\beta,X,\frac{m}{q},\frac{n}{q})\overline{\lambda_{\gamma'}(\beta,X,\frac{m'}{q},\frac{n'}{q})}e((E_\gamma-E_{\gamma'})\beta)d\beta.
\end{align*}
\begin{remark} To estimate $\lambda_{\gamma}(\beta,X,\frac{m}{q},\frac{n}{q})$, we restrict $m,n$ in the modular group such that $-1/2<m/q,n/q\leq 1/2$. By Poisson summation formula,
\begin{equation}\label{lambdagamma}
\begin{aligned}
&\sum_{x\in \mathbb{Z}}\Psi(\frac{x}{X})e\left(A_\gamma \beta x^2+(C_\gamma \beta-\frac{m}{q}) x\right)=\sum_{\xi \in \mathbb{Z}}\int_{-\infty}^{\infty}\Psi(\frac{x}{X})e(A_\gamma \beta x^2+(C_\gamma \beta-\frac{m}{q}) x)e(\xi x)dx\\
=& \sum_{\xi \in \mathbb{Z}}X\int_{-\infty}^{\infty}\Psi(x)e(A_\gamma \beta X^2 x^2+(C_\gamma \beta-\frac{m}{q}) Xx)e(\xi X x)dx.
\end{aligned}
\end{equation}
Noting that $A_\gamma \beta X^2\ll \xi X$, $C_\gamma \beta\ll \xi X$ and $mX/q\ll \xi X$, the integral is $\ll (\xi X)^{-L}$ if $\xi\neq 0$ and $\ll (A_\gamma \beta X^2)^{-1/2}$ if $\xi=0$ by Lemma \ref{NSP} and Lemma \ref{SP}.
The main contribution comes from the term $\xi=0$ and $m<q/X^{1-2\omega}$ ($2\omega$ here is to make sure that $mX/q$ is strictly greater than $1$ with a small power when $m\geq q/X^{1-2\omega}$).
\end{remark}


Split $I_Q$ into $J_Q^=$ and $J_Q^{\neq}$ according to $A_\gamma=A_{\gamma'}$ or not. 
Denote by
\begin{align*}
J_Q^=
=&\sum_{Q<q\leq 2Q}\sum_{\substack{\gamma,\gamma'\\A_\gamma=A_{\gamma'}}}\sum_{m,n(q)}\sum_{m',n'(q)}S(q,\gamma,-m,-n,\gamma',-m',-n')\\
&\times \int_{-1/qM}^{1/qM}\lambda_{\gamma}(\beta,X,\frac{m}{q},\frac{n}{q})\overline{\lambda_{\gamma'}(\beta,X,\frac{m'}{q},\frac{n'}{q})}e((E_\gamma-E_{\gamma'})\beta)d\beta,
\end{align*}
and 
\begin{align*}
J_Q^{\neq}
=&\sum_{Q<q\leq 2Q}\sum_{\substack{\gamma,\gamma'\\A_\gamma\neq A_{\gamma'}}}\sum_{m,n(q)}\sum_{m',n'(q)}S(q,\gamma,-m,-n,\gamma',-m',-n')\\
&\times \int_{-1/qM}^{1/qM}\lambda_{\gamma}(\beta,X,\frac{m}{q},\frac{n}{q})\overline{\lambda_{\gamma'}(\beta,X,\frac{m'}{q},\frac{n'}{q})}e((E_\gamma-E_{\gamma'})\beta)d\beta.
\end{align*}

Firstly we estimate $J_Q^{\neq}$. As before, we write $q=2^\nu q' q''$ with $(2,q'q'')=1$, $(q',q'')=1$, $q'=p_1^{\nu_1}\cdots p_s^{\nu_s}$ with $p_i^{\nu_i}\nmid A_\gamma, A_{\gamma'}$ and $q''\mid {\rm lcm}\{A_\gamma, A_{\gamma'}\}$. Denote $\ell_1=(q'',A_{\gamma}), \ell_2=(q'',A_{\gamma'})$. The bound for $S$ in \eqref{bdSneq}-\eqref{bd2power1} and \eqref{triviallambda} imply
\begin{equation*}
\begin{aligned}
J_Q^{\neq}\ll \frac{X^4}{QM}\sum_{\substack{Q<q<2Q\\q=2^\nu q'q''}}\sum_{\substack{\gamma,\gamma'\\ A_{\gamma}\neq A_{\gamma'}}}\sum_{m,n\leq q/X^{1-2\omega}}\sum_{m',n'\leq q/X^{1-2\omega}}&(q',A_{\gamma})(q',A_{\gamma'}) q'^{-5/4+\ve}(q',A_{\gamma}-A_{\gamma'})^{1/4}\\&\times 2^\nu \cdot {\bf 1} \{2^{\nu-1}\mid (A_{\gamma},A_{\gamma'})\}\frac{\ell_1\ell_2}{q''}.
\end{aligned}
\end{equation*}
Bound $(q',A_{\gamma})$, $(q',A_{\gamma'})$, $(q',A_{\gamma}-A_{\gamma'})$ by $T^2$. Then we obtain
\begin{equation}\label{Jneq}
\begin{aligned}
J_Q^{\neq}
\ll& \frac{X^4}{QM}QT^{4\delta}\frac{Q^4}{X^{4-8\omega}}T^4Q^{-5/4+\ve}T^{15/2}T^{1/2}T^4
\ll T^{4\delta-4}N^{1-\eta_0},
\end{aligned}
\end{equation}
if $X\gg T^{640}$.

We divide $J_Q^=$ into $J_Q^{=,=}$ and $J_Q^{=,\neq}$ according to $\mathfrak{D}m^2+n^2=\mathfrak{D}m'^2+n'^2$ or not. Thanks to \eqref{=neq}-\eqref{bd2power1}, we have
\begin{equation}\label{J=neq}
\begin{aligned}
J_Q^{=,\neq}\ll& \frac{X^4}{QM} \sum_{\substack{\gamma,\gamma'\\ A_{\gamma}=A_{\gamma'}}}\sum_{\substack{Q<q\leq 2Q\\q=2^\nu q'q''}}\sum_{m,n,m',n'\leq q/X^{1-2\omega}}q''^{-1}\ell_1\ell_2 (q',A_{\gamma})^2q'^{-5/4+\ve}\\
&\times(q',(\mathfrak{D}m^2+n^2)-(\mathfrak{D}m'^2+n'^2))^{1/4}\times2^\nu \cdot {\bf 1} \{2^{\nu-1}\mid A_{\gamma}\}\\
\ll& \frac{X^4}{QM}T^{2\delta+2+\ve}Q\frac{Q^4}{X^{4-8\omega}}T^6 Q^{-5/4}(T^6)^{5/4}\frac{Q^{1/2}}{X^{1/2-\omega}}T^2
\ll
 T^{4\delta-4}N^{1-\eta_0},
\end{aligned}
\end{equation}
since we have a power saving for $X$ if $X\gg T^{600-8\delta}$.

For $J_Q^{=,=}$, we can not get enough saving with the bound for $S$ only. More saving  from a more delicate estimate is needed to get the desired bound. 
Write $q=2^\nu q' q''$ as before.
\begin{align*}
J_Q^{=,=}= &\sum_{\substack{\gamma,\gamma'\\ A_{\gamma}=A_{\gamma'}}}\sum_{\substack{Q<q \leq 2Q\\q=2^{\nu}q'q''}}\sum_{m,n(q)}\sum_{\substack{m',n'(q)\\\mathfrak{D}m^2+n^2=\mathfrak{D}m'^2+n'^2}}S(2^{\nu}q'',\gamma,m,n,\gamma',m',n')S(q',\gamma,m,n,\gamma',m',n')\\
&\times \int_{-1/qM}^{1/qM}\lambda_{\gamma}(\beta,X,\frac{m}{q},\frac{n}{q})\overline{\lambda_{\gamma'}(\beta,X,\frac{m'}{q},\frac{n'}{q})}e((E_\gamma-E_{\gamma'})\beta)d\beta.
\end{align*}

Firstly, we deal with the case $Q<T^{2-1/100}X$.
As before, main contribution comes from terms of $m,n, m',n' \leq q/X^{1-2\omega}$. Let $d=(A_\gamma,q')$.
According to the definition of $\chi(\cdot)$ in \eqref{defchi},
consider congruence equations
\begin{equation}\label{congeq}
\begin{aligned}
\begin{cases}
rC_\gamma-m &\equiv 0\quad (\mod d),\\
rD_\gamma-n&\equiv 0\quad  (\mod d),\\
rC_{\gamma'}-m'&\equiv 0\quad (\mod d),\\
rD_{\gamma'}-n'&\equiv 0\quad (\mod d).
\end{cases}
\end{aligned}
\end{equation}
For fixed $r$ and $\gamma$, there is at most one $m$ (resp. $n$) modulo $d$ satisfying the above congruence equations.  The number of $\gamma$ satisfying $2^{\nu-1} q'' d\mid A_{\gamma}$ is $\ll T^4/4^{\nu}q''^2 d^2$, and the number of $\gamma'$ 
with $A_\gamma=A_{\gamma'}$ is $\ll T^{2+\epsilon}$.
By using \eqref{s1s2}-\eqref{s2}, \eqref{trivialbdS}, \eqref{bd2power0} and \eqref{bd2power1}, we have, for $Q<T^{2-1/100}X$,
\begin{equation}\label{J==1}
\begin{aligned}
J_Q^{=,=}\ll &\sum_{d\leq 2T^2}\sum_{\substack{Q<2^{\nu}q'q''\leq 2Q\\ d\mid q'} }\sum_{\substack{\gamma\\ 2^{\nu-1}q'' d\mid A_{\gamma}}}\sum_{\substack{\gamma'\\ A_{\gamma}=A_{\gamma'}}}\frac{Q^2}{X^{2-4\omega}d^2}2^{\nu}q''\frac{d^2}{q'}\frac{X^4}{QM}\\
\ll &\frac{X^2}{M}\sum_{d\leq 2T^2}\frac{QT^{6+\ve}}{d^3} X^{4\omega}
\ll N^{1-\eta_0}T^{4\delta-4},
\end{aligned}
\end{equation}
if $\delta \geq 2-1/400$ ignoring $\omega$ sufficiently small.
(Here for $S(q',\gamma,m,n,\gamma',m',n')$, we use the trivial bound in \eqref{trivialbdS} instead of \eqref{=neq}. In this way, we can reverse the order of sums over $r$ and $m,n,m',n'$ to get the saving by reducing the numbers of $m,n,m',n'$. )

In the following, we shall deal with the case $Q\geq T^{2-1/100}X$. Let $H_0$ be a parameter remaining at our disposition. As before, write $q=2^{\nu}p_1^{\nu_1}\cdots p_s^{\nu_s}\times p_{s+1}^{\nu_{s+1}}\cdots p_{w}^{\nu_{w}}$ where $(p_i,p_j)=1$ for $i\neq j$, $p_i^{\nu_i}\mid A_{\gamma}$ for $1\leq i\leq s$ and $p_i^{\nu_i}\nmid A_{\gamma}$ for $s+1\leq i\leq w$. Denote by
$q_1=2^\nu$, $q_2=p_1^{\nu_1}\cdots p_s^{\nu_s}$ and $q_3=p_{s+1}^{\nu_{s+1}}\cdots p_{w}^{\nu_{w}}$.

Denote $d=\gcd(A_\gamma, q_3)$. Divide $J_Q^{=,=}$ into two parts according to $dq_1q_2>H_0$ or not.
Denote by
\begin{align*}
\mathscr{D}_1= &\sum_{d\leq 2T^2}\sum_{\substack{\gamma\\ d\mid A_{\gamma}}}\sum_{\substack{\gamma'\\ A_{\gamma}=A_{\gamma'}}}\sum_{q_1}\sum_{q_2\geq H_0/dq_1}\sum_{\substack{d\mid q_3\\ Q<q_1q_2q_3\leq 2Q}}\sum_{\substack{m,n,m',n'\leq q_1q_2q_3/X^{1-2\omega}\\ \mathfrak{D}m^2+n^2=\mathfrak{D}m'^2+n'^2}}\!\!\!\!\!\!\!\!|S(q_1q_2q_3,\gamma,m,n,\gamma',m',n')|\\
&\times \left|\int_{-1/q_1q_2q_3M}^{1/q_1q_2q_3M}\lambda_{\gamma}(\beta,X,\frac{m}{q_1q_2q_3},\frac{n}{q_1q_2q_3})\overline{\lambda_{\gamma'}(\beta,X,\frac{m'}{q_1q_2q_3},\frac{n'}{q_1q_2q_3})}e((E_\gamma-E_{\gamma'})\beta)d\beta\right|.
\end{align*}

According to Corollary \ref{corSodd} and Corollary \ref{corSeven}, $q_1\mid 2A_\gamma$, $q_2\mid A_\gamma$, $S(q_1,\gamma,m,n,\gamma',m',n')\ll q_1$ and $S(q_2,\gamma,m,n,\gamma',m',n')\ll q_2$.
Similar argument as above (modulo $dq_1q_2/2$ in congruence equations \eqref{congeq} here), we have
\begin{equation}\label{mD1}
\begin{aligned}
\mathscr{D}_1&\ll \sum_{H_0<dq_1q_2\leq 2T^2}\frac{T^{6+\ve}}{(dq_1q_2)^2}\sum_{d\mid q_3, Q<q_1q_2q_3\leq 2Q}\frac{(q_1q_2q_3)^2}{d^2q_1^2q_2^2X^{{2-4\omega}}}q_1q_2\frac{d^2}{q_3}\frac{X^4}{QM}\\
&\ll \sum_{H_0<dq_1q_2\leq 2T^2}\frac{T^6}{(dq_1q_2)^2}\sum_{q_3}\frac{(q_1q_2q_3)^2}{d^2q_1^2q_2^2X^{{2-4\omega}}}q_1q_2\frac{1}{q_3}d^2\frac{X^4}{QM}\\
&\ll NT^{4\delta-4}\frac{T^{8-4\delta}}{H^2_0}X^{4\omega}\ll N^{1-\eta_0}T^{4\delta-4},
\end{aligned}
\end{equation}
if $H_0 \geq T^{4-2\delta}.$

It remains to deal with
\begin{align*}
\mathscr{D}_2= &\sum_{d\leq 2T^2}\sum_{\substack{\gamma\\ d\mid A_{\gamma}}}\sum_{\substack{\gamma'\\ A_{\gamma}=A_{\gamma'}}}\sum_{dq_1q_2\leq H_0}\sum_{\substack{(A_\gamma,q_3)=d\\Q<q_1q_2q_3\leq 2Q}}\sum_{\substack{m,n,m',n'\leq q_1q_2q_3/X^{1-2\omega}\\ \mathfrak{D}m^2+n^2=\mathfrak{D}m'^2+n'^2}}S(q_1q_2q_3,\gamma,m,n,\gamma',m',n')\\
&\times \int_{-1/q_1q_2q_3M}^{1/q_1q_2q_3M}\lambda_{\gamma}(\beta,X,\frac{m}{q_1q_2q_3},\frac{n}{q_1q_2q_3})\overline{\lambda_{\gamma'}(\beta,X,\frac{m'}{q_1q_2q_3},\frac{n'}{q_1q_2q_3})}e((E_\gamma-E_{\gamma'})\beta)d\beta.
\end{align*}
Rewriting the summation, we can get
\begin{align*}
\mathscr{D}_2= &\sum_{\substack{\gamma,\gamma'\\ A_{\gamma}=A_{\gamma'}}}\sum_{\substack{m,n,m',n'\leq 2T^2X^{2\omega}\\ \mathfrak{D}m^2+n^2=\mathfrak{D}m'^2+n'^2}}\sum_{\substack{d \mid A_{\gamma}\\dq_1q_2 \leq H_0}}\sum_{\substack{Q<q_1q_2q_3\leq 2Q,(q_3,A_\gamma)=d\\ q_1q_2q_3/X^{1-2\omega}>\max\{m,n,m',n'\}}}S(q_1q_2q_3,\gamma,m,n,\gamma',m',n')\\
&\times \int_{-1/q_1q_2q_3M}^{1/q_1q_2q_3M}\lambda_{\gamma}(\beta,X,\frac{m}{q_1q_2q_3},\frac{n}{q_1q_2q_3})\overline{\lambda_{\gamma'}(\beta,X,\frac{m'}{q_1q_2q_3},\frac{n'}{q_1q_2q_3})}e((E_\gamma-E_{\gamma'})\beta)d\beta.
\end{align*}
%

In order to divide the inner sum of $q_3$ into congruence classes modulo $A_\gamma$, we shall substitute the integral by another one independent of $q_3$ with a negligible error.  

Notice that for any $0<t<A_\gamma^2$, by the definition of $\lambda_{\gamma}(\beta,X,s,t)$,
\begin{align*}
\lambda_{\gamma}(\beta,X,\frac{m}{q},\frac{n}{q})-\lambda_{\gamma}(\beta,X,\frac{m}{q+t},\frac{n}{q+t})\ll \frac{\max\{m,n\}X^3t}{Q^2}\ll X^{1+2\omega}T^{2+1/50},
\end{align*}
by using the fact
$$e\left(-\frac{m}{q+t}x\right)=e\left(-\frac{m}{q}x\right)+O\left(\frac{mtx}{q(q+t)}\right).$$
For fixed $\gamma, \gamma', m,n,m',n'$, let $\max\{m,n,m',n',\frac{Q}{X^{1-2\omega}}\}X^{1-2\omega}\leq Q'\leq 2Q$.  We have, for $Q'\leq q < Q'+A_\gamma^2$,
\begin{align*}
& \int_{-1/qM}^{1/qM}\lambda_{\gamma}(\beta,X,\frac{m}{q},\frac{n}{q})\overline{\lambda_{\gamma'}(\beta,X,\frac{m'}{q},\frac{n'}{q})}e((E_\gamma-E_{\gamma'})\beta)d\beta
\\ 
=& \int_{-1/qM}^{1/qM}\lambda_{\gamma}(\beta,X,\frac{m}{Q'},\frac{n}{Q'})\overline{\lambda_{\gamma'}(\beta,X,\frac{m'}{Q'},\frac{n'}{Q'})}e((E_\gamma-E_{\gamma'})\beta)d\beta+O\left(\frac{X^{3+2\omega}T^{2+1/50}}{QM}\right)\\
=& \int_{-1/Q'M}^{1/Q'M}\lambda_{\gamma}(\beta,X,\frac{m}{Q'},\frac{n}{Q'})\overline{\lambda_{\gamma'}(\beta,X,\frac{m'}{Q'},\frac{n'}{Q'})}e((E_\gamma-E_{\gamma'})\beta)d\beta+O\left(\frac{X^{3+2\omega}T^{2+1/50}}{QM}+\frac{X^4T^4}{Q^2M}\right). 
\end{align*}
Then we get (letting $Q'=Q+2jq_1q_2 A_\gamma$ in the above)
\begin{equation}\label{I==}
\begin{aligned}
\mathscr{D}_2= &\sum_{\substack{\gamma,\gamma'\\ A_{\gamma}=A_{\gamma'}}}\sum_{\substack{m,n,m',n'\leq 2T^2X^{2\omega}\\ \mathfrak{D}m^2+n^2=\mathfrak{D}m'^2+n'^2}}\sum_{\substack{d \mid A_{\gamma} \\ dq_1q_2\leq H_0}}{\sum_{j}}^*\sum_{\substack{\frac{Q}{q_1q_2}+2jA_{\gamma}<q_3\leq \frac{Q}{q_1q_2}+2(j+1)A_{\gamma} \\ (q_3,A_\gamma)=d}}\!\!\!\!\!\!\!S(q_1q_2q_3,\gamma,m,n,\gamma',m',n')\\
&\times \int_{-1/Q'M}^{1/Q'M}\lambda_{\gamma}(\beta,X,\frac{m}{Q'},\frac{n}{Q'})\overline{\lambda_{\gamma'}(\beta,X,\frac{m'}{Q'},\frac{n'}{Q'})}e((E_\gamma-E_{\gamma'})\beta)d\beta+O(\fE_1),
\end{aligned}
\end{equation}
where the sum $\sum^*$ runs over all $j\geq 0$ satisfying
$$\frac{Q}{q_1q_2}+2jA_\gamma>X^{1-2\omega}\max\{m,n,m',n'\}/q_1q_2, \qquad 2(j+1)A_\gamma\leq \frac{2Q}{q_1q_2}.$$
Here (by using $S(q_1q_2q_3,\gamma,m,n,\gamma',m',n')\ll q_1q_2d^2/q_3$)
\begin{align*}
\fE_1\ll &\sum_{\substack{\gamma,\gamma'\\ A_{\gamma}=A_{\gamma'}}}\sum_{\substack{m,n,m',n'\leq 2T^2X^{2\omega}\\ \mathfrak{D}m^2+n^2=\mathfrak{D}m'^2+n'^2}}\sum_{\substack{d \mid A_{\gamma}\\q_1q_2d\leq H_0}}\sum_{\substack{Q<q_1q_2q_3\leq 2Q,(q_3,A_\gamma)=d\\ q_1q_2q_3/X^{1-2\omega}>\max\{m,n,m',n'\}}}\frac{q_1q_2d^2}{q_3}\\
&\times \left(\frac{X^{3+2\omega}T^{2+1/50}}{QM}+\frac{X^4T^4}{Q^2M}\right)+ \sum_{\substack{\gamma,\gamma'\\ A_{\gamma}=A_{\gamma'}}}\sum_{\substack{m,n,m',n'\leq 2T^2X^{2\omega}\\ \mathfrak{D}m^2+n^2=\mathfrak{D}m'^2+n'^2}}\sum_{d \mid A_{\gamma}, q_1q_2d\leq H_0}\frac{q_1q_2d^2}{q_3}\frac{X^4}{QM}T^2.
\end{align*}
It is easy to see that we have saving in $X$. Then
$$\fE_1 \ll N^{1-\eta_0}T^{4\delta-4}.$$
%

Write $q_3=dq'$. Let $Q'=Q+2jq_1q_2A_\gamma$ and $Q''=\frac{Q}{q_1q_2}+2jA_\gamma$. Denote by
\begin{align*}
\mathscr{S}(\gamma,m,n,\gamma',m',n'):=\sum_{d\mid A_\gamma,q_1q_2d\leq H_0}{\sum_j}^*\tilde{I}(j,q_1,q_2)\sum_{\substack{Q''<dq'\leq Q''+2A_\gamma\\ (q', A_\gamma/d)=1}}S(q_1q_2q_3,\gamma,m,n,\gamma',m',n'),
\end{align*}
with
$$\tilde{I}(j,q_1,q_2)= \int_{-1/Q'M}^{1/Q'M}\lambda_{\gamma}(\beta,X,\frac{m}{Q'},\frac{n}{Q'})\overline{\lambda_{\gamma'}(\beta,X,\frac{m'}{Q'},\frac{n'}{Q'})}e((E_\gamma-E_{\gamma'})\beta)d\beta.$$
According to \eqref{s1s2}, \eqref{s1} and \eqref{s2}, we can write
\begin{align*}
\mathscr{S}(\gamma,m,n,\gamma',m',n')=\!\!\!\sum_{d\mid A_\gamma,q_1q_2d\leq H_0}\!\!\!\!\!S(q_1q_2,\gamma,m,n,\gamma',m',n'){\sum_j}^*\tilde{I}(j,q_1,q_2)\!\!\!\sum_{\substack{Q''/d<q'\leq (Q''+2A_\gamma)/d\\ (q', A_\gamma/d)=1}}S_1S_2.
\end{align*}
It should be mentioned that, we cannot take the absolute value of $S$ and $\tilde{I}$ as before, since we are going to get saving from $S_1$. For $S_2$, we can consider the congruence equation system in \eqref{defchi} (see also \eqref{congeq}) again. 
({It seems we can get nothing from $\mathfrak{D}m^2+n^2=\mathfrak{D}m'^2+n'^2$. In fact, $r^2(\mathfrak{D}C^2_{\gamma'}+D^2_{\gamma'}-\mathfrak{D}C^2_{\gamma}-D^2_{\gamma})\equiv 0 \ (\mod d)$ since $\mathfrak{D}C^2_{\gamma}+D^2_{\gamma}=4\mathfrak{D}(d_1^2+\mathfrak{D}d_2^2)A_{\gamma}$.}) If $d=p_1^{\nu_1}\cdot p_s^{\nu_s}$, denote
$d'=p_1^{\nu'_1}\cdot p_s^{\nu'_s}$ where $\nu'_i=\nu_i$ if $p\nmid \mathfrak{D}$ and $\nu'_j=\nu_j+1$ if $p \mid \mathfrak{D}$. If fact, $d'=(q_3, B_\gamma)$. Then $d\mid d'$ and $\frac{d'}{d}\mid \mathfrak{D}$.
Let $v_1=\gcd(d, C_\gamma)$, $v_2=\gcd(d', D_\gamma)$, $v_3=\gcd(d, C_{\gamma'})$, $v_4=\gcd(d', D_{\gamma'})$.
According to the definition of $S_2$, by Chinese Remaining Theorem, for fixed $C_\gamma, D_\gamma, C_{\gamma'}, D_{\gamma'}$ and $m, n, m', n', d,d'$, there is one and only one $r$ (appearing when $C_\gamma, D_\gamma, C_{\gamma'}, D_{\gamma'}$ and $m, n, m', n', d,d'$ satisfy certain conditions) modulo 
$$D=[\frac{d}{v_1},\frac{d'}{v_2},\frac{d}{v_3},\frac{d'}{v_4}],$$
satisfying $\chi(\cdot)=1$. Thus, $|S_2|\ll q_3/D=dq'/D$.
Observe that $S_2$ is a function of $C_\gamma, D_\gamma, C_{\gamma'}, D_{\gamma'}$ and $m, n, m', n', d,d'$ and independent of $q'$. By \eqref{s1}, we obtain
\begin{equation}\label{eq1}
\begin{aligned}
&\mathscr{S}(\gamma,m,n,\gamma',m',n')
\ll&\sum_{d\mid A_\gamma,q_1q_2d\leq H_0}\frac{d'}{D}|S(q_1q_2,\gamma,m,n,\gamma',m',n')|{\sum_j}^*|\tilde{I}(j,q_1,q_2)||\mathscr{S}'|,
\end{aligned}
\end{equation}
with
$$\mathscr{S}':=\!\!\!\!\!\!\sum_{\substack{Q''/d<q'\leq Q''/d+2A_\gamma/d\\ (q', A_\gamma/d)=1}}\frac{1}{q'}e\left(-(\mathfrak{D}C_\gamma m+D_\gamma n)/d'\cdot\frac{\overline{2A_\gamma/d} }{q_3}\right)e\left((\mathfrak{D}C_{\gamma'} m'+D_{\gamma'}n')/d'\cdot \frac{\overline{2A_{\gamma'}/d}}{q_3}\right).$$
Observe that for some $k$
$$2A_\gamma/d \cdot \overline{2A_\gamma/d}=kq_3+1,$$
then
$$\frac{\overline{2A_\gamma/d}}{q_3}=\frac{k}{2A_\gamma/d}+\frac{1}{q_3\cdot 2A_\gamma/d}.$$
For fixed $d$ and a given $q'\, (\mod 2A_\gamma/d)$, we have a unique $k\, (\mod 2A_\gamma/d)$ satisfying the above equation. Therefore, by using the fact
$$\frac{1}{q'}-\frac{1}{Q''/d}\ll \frac{A_\gamma d}{Q''^2},$$
we can write
\begin{align*}
\mathscr{S}'&={\sum_{k(\mod 2A_\gamma/d)}}'\frac{1}{q'}e\left(-\frac{kG}{2A_\gamma/d}\right)\left(1+O\left(\frac{dG}{Q''A_\gamma}\right)\right)\\
&= \frac{1}{Q''/d}{\sum_{k(\mod 2A_\gamma/d)}}'e\left(-\frac{kG}{2A_\gamma/d}\right)+O\left(\frac{dG+A_\gamma^2}{Q''^2}\right),
\end{align*}
with
$$G=(C_{\gamma} m+D_{\gamma}n-C_{\gamma'} m'-D_{\gamma'}n')/d',$$
and a reduced modulo group in the last sum.
Let $\ell_0=\gcd(2A_\gamma/d, G)$.  
By using \eqref{ramanujansumbd},
we have
\begin{align*}
\mathscr{S}'&\ll \frac{\ell_0}{Q''/d}.
\end{align*}

Inserting the above estimates into \eqref{eq1} and then into \eqref{I==}, one obtains
\begin{equation*}
\mathscr{D}_2\ll \frac{X^4}{QM}\sum_{\substack{\gamma,\gamma'\\ A_{\gamma}=A_{\gamma'}}}\sum_{m,n\leq 2T^2X^{2\omega}}\sum_{\substack{m',n'\leq 2T^2X^{2\omega}\\ \mathfrak{D}m^2+n^2=\mathfrak{D}m'^2+n'^2}}\sum_{\substack{d\mid A_{\gamma}\\dq_1q_2\leq H_0}}{\sum_{j}}^*\frac{d'}{D}q_1q_2\frac{\ell_0}{Q''/d}.
\end{equation*}
By using $D\geq 1$ and Lemma \ref{samebottom}, we can bound $\mathscr{D}_2$ by
\begin{align*}
\mathscr{D}_2\ll& \frac{X^4}{QM}\sum_{m,n\leq 2T^2X^{2\omega}}\sum_{\substack{m',n'\leq 2T^2X^{2\omega}\\ \mathfrak{D}m^2+n^2=\mathfrak{D}m'^2+n'^2}}\sum_{\substack{\gamma,\gamma'\\ A_{\gamma}=A_{\gamma'}}}\sum_{\substack{d\mid A_{\gamma},dq_1q_2\leq H_0}}\frac{Q}{q_1q_2A_\gamma}\frac{d'}{D}\ell_0\frac{q_1^2q_2^2}{Q/d}\\
\ll& \frac{X^4}{QM}\sum_{dq_1q_2\leq H_0}dq_1q_2\sum_{\substack{c_1,c_2\leq T\\d\mid c_1^2+c_2^2\mathfrak{D}}}\frac{d}{c_1^2+c_2^2\mathfrak{D}}\sum_{\substack{c_1',c_2'\leq T\\c_1'^2+\mathfrak{D}c_2'^2= c_1^2+\mathfrak{D}c_2^2}}\sum_{\substack{m,n,m',n' \leq 2T^2X^{2\omega}\\ \mathfrak{D}m^2+n^2=\mathfrak{D}m'^2+n'^2}}
\sum_{d_1,d_2,d'_1,d'_2\leq T}\ell_0.
\end{align*}

%

According to \eqref{defQ}, rewrite $G$ as the form
$$G=[(2c_1m-2\mathfrak{D}c_2n)d_1+(2c_2m+2c_1n)\mathfrak{D}d_2-(2c'_1m'-2\mathfrak{D}c'_2n')d'_1-(2c'_2m'+2c'_1n')\mathfrak{D}d'_2]/d'.$$
If we fix $c_1,c_2,m,n$ and $c'_1,c'_2,m',n',d_1',d_2'$, the number of solutions of $(d_1,d_2)$ for the congruence equation
$$(2c_1m-2\mathfrak{D}c_2n)d_1+(2c_2m+2c_1n)\mathfrak{D}d_2 \equiv (2c'_1m'-2\mathfrak{D}c'_2n')d'_1+(2c'_2m'+2c'_1n')\mathfrak{D}d'_2 \quad (\mod \ell_0)$$
is
$$\ll \frac{T^2}{\ell_0/(2c_1m-2\mathfrak{D}c_2n,(2c_2m+2c_1n)\mathfrak{D},\ell_0)}+T,$$
thanks to Lemma \ref{solcong}.

Thus, writing $\ell_1=(2c_1m-2\mathfrak{D}c_2n,(2c_2m+2c_1n)\mathfrak{D},\ell_0)$, we have
\begin{align*}
\mathscr{D}_2\ll&  \frac{X^4}{QM}H_0\sum_{dq_1q_2\leq H_0}\sum_{\substack{c_1,c_2,c'_1,c'_2\leq T\\d\mid c_1^2+c_2^2\mathfrak{D}\\c_1^2+c_2^2\mathfrak{D}=c_1'^2+c_2'^2\mathfrak{D}}}\frac{d}{c_1^2+c_2^2\mathfrak{D}}\sum_{\substack{\ell_0\mid c_1^2+c_2^2\mathfrak{D}}}\sum_{\substack{m,n,m',n' \leq 2T^2X^{2\omega}\\ \mathfrak{D}m^2+n^2=\mathfrak{D}m'^2+n'^2}}\sum_{d'_1,d'_2\leq T}T^2\ell_1\\
&+\frac{X^4}{QM}H_0\sum_{dq_1q_2\leq H_0}\sum_{\substack{c_1,c_2,c'_1,c'_2\leq T\\d\mid c_1^2+c_2^2\mathfrak{D}\\c_1^2+c_2^2\mathfrak{D}=c_1'^2+c_2'^2\mathfrak{D}}}\frac{d}{c_1^2+c_2^2\mathfrak{D}}\sum_{\substack{\ell_0\mid c_1^2+c_2^2\mathfrak{D}}}\sum_{\substack{m,n,m',n' \leq 2T^2X^{2\omega}\\ \mathfrak{D}m^2+n^2=\mathfrak{D}m'^2+n'^2}}\sum_{d'_1,d'_2\leq T}T\ell_0\\
=:&\mathscr{D}_{21}+\mathscr{D}_{22}.
\end{align*}
The second sum is
\begin{equation}\label{D22}
\begin{aligned}
\mathscr{D}_{22}\ll& \frac{X^4}{QM}H_0\sum_{dq_1q_2\leq H_0}\sum_{\ell_0\leq T^2}\sum_{n\ll T^2,d\ell_0\mid n}\frac{d\ell_0}{n}\sum_{\substack{c_1,c_2,c'_1,c'_2\leq T\\n=c_1^2+c_2^2\mathfrak{D}=c_1'^2+c_2'^2\mathfrak{D}}}\sum_{\substack{m,n,m',n' \leq 2T^2X^{2\omega}\\ \mathfrak{D}m^2+n^2=\mathfrak{D}m'^2+n'^2}}T^3\\
\ll& \frac{X^4}{QM}H^2_0T^2T^4X^{4\omega}T^3 \ll NT^{3+1/100}H_0^2\ll T^{4\delta-4}N^{1-\eta_0},
\end{aligned}
\end{equation}
if $H_0<T^{2\delta-7/2-1/200}$.

For the first sum, we have
\begin{align*}
\mathscr{D}_{21}\ll&  \frac{X^4}{QM}H_0\sum_{dq_1q_2\leq H_0}\sum_{\substack{c_1,c_2,c'_1,c'_2\leq T\\d\mid c_1^2+c_2^2\mathfrak{D},c_1^2+c_2^2\mathfrak{D}=c_1'^2+c_2'^2\mathfrak{D}}}\frac{d}{c_1^2+c_2^2\mathfrak{D}}\sum_{\substack{\ell_0\mid c_1^2+c_2^2\mathfrak{D}}}\sum_{\ell_1\mid \ell_0}T^4\ell_1\\
&\times \sum_{m,n\leq 2T^2X^{2\omega}}1\{\ell_1\mid 2c_1m-2\mathfrak{D}c_2n\}\sum_{\substack{m',n'\\m^2+n^2=m'^2+n'^2}}1.
\end{align*}
Using Lemma \ref{solcong} again, we obtain
\begin{align*}
\mathscr{D}_{21}\ll&  \frac{X^4T^{4+\ve}}{QM}H_0\sum_{dq_1q_2\leq H_0}\sum_{\substack{c_1,c_2,c'_1,c'_2\leq T\\d\mid c_1^2+c_2^2\mathfrak{D},c_1^2+c_2^2\mathfrak{D}=c_1'^2+c_2'^2\mathfrak{D}}}\frac{d}{c_1^2+c_2^2\mathfrak{D}}\sum_{\substack{\ell_0\mid c_1^2+c_2^2\mathfrak{D} }}\sum_{\ell_1\mid \ell_0}\frac{\ell_1 \cdot T^4X^{4\omega}}{\ell_1/(c_1,c_2\mathfrak{D},\ell_1)}\\
\ll&\frac{X^4T^{8+\ve}}{QM}X^{4\omega}H_0\sum_{dq_1q_2\leq H_0}\sum_{\substack{c_1,c_2,c'_1,c'_2\leq T\\d\mid c_1^2+c_2^2\mathfrak{D},c_1^2+c_2^2\mathfrak{D}=c_1'^2+c_2'^2\mathfrak{D}}}\frac{d}{c_1^2+c_2^2\mathfrak{D}}\sum_{\substack{\ell_0\mid c_1^2+c_2^2\mathfrak{D}}}\sum_{\ell_1\mid \ell_0}(c_1,c_2),
\end{align*}
which implies that
\begin{align*}
\mathscr{D}_{21}
\ll \frac{X^4T^{8+\ve}}{QM}X^{4\omega}H_0\sum_{dq_1q_2\leq H_0}\sum_{\substack{c_1,c_2\leq T\\d\mid c_1^2+c_2^2\mathfrak{D}}}\frac{d}{c_1^2+c_2^2\mathfrak{D}}(c_1,c_2).
\end{align*}
Let $t=(c_1,c_2)$. Then
\begin{align*}
\mathscr{D}_{21} \ll& \frac{X^4T^{8+\ve}}{QM}X^{4\omega}H_0\sum_{dq_1q_2\leq H_0}\sum_{t\leq T}\sum_{c_1,c_2\leq T/t}\frac{d}{t^2(c_1^2+c_2^2\mathfrak{D})}t\\
\ll& \frac{X^4T^{8+\ve}}{QM}X^{4\omega}H_0\sum_{dq_1q_2\leq H_0}d\sum_{t\leq T}\sum_{c_1,c_2\leq T/t}\frac{1}{tc_1c_2}.
\end{align*}
It follows that
\begin{align}\label{mA2}
\mathscr{D}_{21}\ll \frac{X^4T^{8+\ve}}{QM}H_0^{3}X^{4\omega}\ll NT^{2+1/100}H_0^{3}X^{4\omega}<N^{1-\eta_0}T^{4\delta-4},
\end{align}
if $H_0<T^{4\delta/3-2-1/300}$.


By \eqref{Jneq}, \eqref{J=neq}, \eqref{J==1}, \eqref{mD1}, \eqref{D22} and \eqref{mA2}, taking $H_0=T^{1/4}$, we can conclude that
\begin{equation}\label{bdI4}
I_4\ll N^{1-\eta_0}T^{4\delta-4}.
\end{equation}


\section{Proof of the Main Theorem}
Firstly, by \eqref{bdI1}, \eqref{bdI2}, \eqref{bdI3} and \eqref{bdI4}, taking
\begin{equation}\label{KQ}
K_0=Q_0^3, \quad Q_0=\min\left\{T^{\tfrac{2(\delta-\Theta_2)}{3}},T^{\tfrac{\delta-\Theta_1}{70}}\right\},
\end{equation}
we can claim Proposition \ref{minor arcs}.

Let $E(N)$ be the set of admissible integers in $[N,2N]$ but not in $S$, then for each $n \in E(N)$, we have
$$R(n)=0$$
and 
$$\mathcal{E}_N(n)=\mathcal{M}_N(n)\gg \frac{T^{2\delta-2}}{\log n}$$
by Proposition \ref{main arcs}.

Thus we obtain,
\begin{equation}\label{p1}
\sum_{n\in E(N)}|\mathcal{E}_N(n)|^2\gg 
\#E(N)\frac{T^{4\delta-4}}{\log^2 n}.
\end{equation}

On the other hand, we have
\begin{equation}\label{p2}
	\sum_{n\in E(N)}|\mathcal{E}_N(n)|^2\ll \sum_{n\in \mathbb{Z}}|\mathcal{E}_N(n)|^2 \ll T^{4\delta-2}N^{1-\eta_0},
\end{equation}
according to Proposition \ref{minor arcs}.

Combining \eqref{p1} and \eqref{p2}, we conclude that
$$\#E(N)\ll N^{1-\eta_0},$$
ignoring the log factors.

\end{document}